\documentclass[11pt]{article}
\usepackage{amsmath,amssymb,amsfonts} % Typical maths resource packages
\usepackage{graphics}                 % Packages to allow inclusion of graphics
\usepackage{makeidx}
\usepackage{graphicx}
\usepackage{sidecap}
\usepackage{appendix}
\usepackage[caption=false,font=normalsize,labelfont=sf,textfont=sf]{subfig}
\usepackage{algorithmic}
\usepackage[algoruled,vlined]{algorithm2e}
\usepackage{fullpage}
\usepackage{float}

\numberwithin{equation}{section}
\numberwithin{figure}{section}
\numberwithin{table}{section}

\newtheorem{theorem}{Theorem}[section]
\newtheorem{proposition}[theorem]{Proposition}
\newtheorem{corollary}[theorem]{Corollary}

\newtheorem{remark}[theorem]{Remark}
\newtheorem{definition}[theorem]{Definition}
\def \srr{\stackrel{\mathrm{def}}{=}}
\newcommand{\rf}[1]{(\ref{#1})}
\def \ehh{Equation \rf}
\def \eh{Eq. \rf}
\def \z{\mathbb{Z}}
\def \br{\begin{remark}}
\def \er{\end{remark}}
\def \bdd{\begin{definition} }
\def \bpp{\begin{proposition} }
\def \bcc{\begin{corollary} }

\def \edd{\end{definition}}
\def \epp{\end{proposition}}
\def \ecc{\end{corollary} }
\def\proof{\goodbreak\medskip\noindent{\small\bf Proof: }}
\def\eop{{\vrule height7pt width7pt depth0pt}\par\bigskip}

\begin{document}
\title{Directional  wavelet packets originating from  polynomial  splines}

\author{Amir Averbuch (1), Pekka Neittaanm\"aki (2) and Valery  Zheludev (1)  \\
((1) School of Computer Science\\
Tel Aviv University, Tel Aviv 69978, Israel\\~\\
(2) Faculty of Mathematical Information Technology\\
 University of Jyv\"askyl\"a, Finland)}
 \date{ }
\maketitle
\begin{abstract}
The paper presents a versatile library of  quasi-analytic  complex-valued wavelet packets (WPs) which originate from polynomial  splines of arbitrary orders. The real parts of the quasi-analytic  WPs are the regular spline-based orthonormal  WPs designed in \cite{ANZ_book3}. The imaginary parts are the so-called complementary orthonormal  WPs that are derived from the Hilbert transforms of the regular WPs and, unlike the symmetric regular  WPs,  are  antisymmetric. Tensor products of 1D quasi-analytic  WPs  provide a diversity of 2D WPs oriented in multiple directions. For example, a set of the fourth-level WPs comprises 62 different  directions. The properties of the presented WPs are  refined frequency resolution,  directionality of waveforms with unlimited number of orientations,  (anti-)symmetry of waveforms and windowed oscillating structure of waveforms  with a variety of frequencies. Directional WPs have a strong potential to be used in various image processing applications such as restoration of degraded images and extraction of  characteristic features from the images.\end{abstract}

\section{Introduction}\label{sec:s1}
Multimedia images as well as biomedical, seismic, and hyper-spectral images, to name a few, comprise smooth regions, edges oriented in various directions and texture, which can have an oscillating structure. One of the main goals of image processing  is  to reconstruct an image from a degraded data that occurs from example from missing a number of pixels, noise and blurring.  Another goal is to extract a limited number of characteristic features from the image for pattern recognition and machine learning applications. Achieving the above goals relies on the fact that practically all images to be processed have a sparse representation in a proper transformed domain. The sparse representation of an image means that it can be approximated by a linear combination of a relatively small number of 2D ``basic" elements (caaled dictionary), while retaining the above mentioned components of the image. The dictionary of such elements should comprise waveforms that \textbf{1. Are oriented in multiple directions }(for capturing edges),\textbf{ 2. Have oscillating structure with multiple frequencies }for retaining texture patterns, and \textbf{3. Have vanishing moments, at least locally }for sparse representation of smooth regions. In addition, properties of the waveform such as
\textbf{4. Refined frequency  separation}, and \textbf{5. Good localization in the spatial domain} are desirable. Last but not least is  \textbf{6. Fast implementation of the corresponding transform s}.

In recent years a number of dictionaries elements  that meet some of the above requirements are constructed in the literature and used in image processing such as pseudo polar \cite{averbuch2008frameworkI,averbuch2008frameworkII}, contourlets \cite{Contour}, curvelets \cite{curve,curve1} and shearlets \cite{kuty,shear}. These dictionaries aare used in various image processing  applications
%One of the most  successful applications of the shearlet transforms to  the image processing  is the methodology based on the  \emph{Digital
such as Affine Shear transforms (DAS-1) \cite{zhuang}. However, while these successfully capture edges in images, these  dictionaries did not demonstrate a satisfactory  texture restoration due to  the lack of oscillating waveforms in the dictionaries.

Another approach to the design of directional  dictionaries consists of the tensor multiplication of complex wavelets (\cite{king1,barakin}), wavelet  frames and wavelet packets (WPs) \cite{jalob1,bay_sele,bhan_zhao,bhan_com_sup, bhan_zhao_zhu}, to name a few.  The tight tensor-product complex  wavelet  frames (TP\_$\mathbb{C}$TF$_{n}$) with different  number of  directions, are designed in \cite{bhan_zhao,bhan_zhao_zhu,bhan_com_sup} and some of them, in particular cptTP\_$\mathbb{C}$TF$_{6}$ TP\_$\mathbb{C}$TF$_{6}$ and TP\_$\mathbb{C}$TF$^{\downarrow}_{6}$, demonstrate impressive performance for image denoising and impainting.
The waveforms in these frames are oriented in 14 directions and, due to the 2-layer structure of their spectra, they possess some oscillatory
 properties.

Some of  disadvantages of the above 2D  TP\_$\mathbb{C}$TF$_{6}$ and TP\_$\mathbb{C}$TF$^{\downarrow}_{6}$ frames such as, for example, limited and fixed number of  directions (14  directions at each decomposition  level) are overcome in
\cite{che_zhuang} (algorithm   \emph{Digital Affine Shear Filter Transform with 2-Layer Structure (DAS-2)}) by the incorporation of the two-layer structure, which is inherent in the TP\_$\mathbb{C}$TF$_{6}$  frames, into directional filter banks introduced in \cite{zhuang}. This improves the  performance of DAS-2 compared to TP\_$\mathbb{C}$TF$_{6}$ on texture-rich images such as ``Barbara", which is not the case for smoother images like ``Lena".

Our motivation was to design a family of dictionaries elements  that maximally meet the requirements 1--6 to utilize them in image processing  applications. For such a design, we have two libraries of orthonormal  WPs originating from the discrete  and the so-called discrete-time splines\footnote{The discrete-time splines are derived by the discretization of polynomial  splines.} of multiple orders (see \cite{ANZ_book3}). The waveforms in both libraries are symmetric, well localized in time domain, their shapes vary from low-frequency  smooth curves to high-frequency  oscillating transients. They can have any number of local vanishing moments (to be defined in Section \ref{sec:ss24}).  Their spectra provide a variety of refined splits of the frequency  domain and shapes of the magnitude spectra tend to a rectangular as the spline's order increases.
Their tensor products possess similar properties extended to 2D setting while disadvantageously, lack directionality.

The following steps are used to design directional WPs: 1. Apply the Hilbert transform  (HT) to the set $\left\{\psi\right\}$ of  orthonormal  WPs thus producing the set  $\left\{\theta=H(\psi)\right\}$.
2. A slight correction of lowest- and highest-frequency   waveforms from the  set  $\left\{\theta\right\}$ provides an orthonormal  set   $\left\{\varphi\right\}$ of the so-called complimentary WPs (cWPs), which are anti-symmetric and whose magnitude spectra coincides with the magnitude spectra of respective WPs from the  set $\left\{\psi\right\}$. 3. Define two sets of complex quasi-analytic  WPs (qWPs) $\left\{\Psi_{+}\srr\psi+i\,\varphi\right\}$ and $\left\{\Psi_{-}\srr\psi-i\,\varphi\right\}$ whose spectra are localized in the positive and negative half-bands of the frequency  domain, respectively. 4. Define two sets of 2D complex qWPs by the tensor multiplication of the qWPs $\left\{\Psi_{\pm}\right\}$ as: $\left\{\Psi_{++}\srr\Psi_{+}\bigotimes\Psi_{+}\right\}$ and $\left\{\Psi_{+-}\srr\Psi_{+}\bigotimes\Psi_{-}\right\}$. 5. The dictionaries we are looking for are obtained as real parts of these qWPs: $\left\{\vartheta_{+}\srr\mathfrak{Re}(\Psi_{++})\right\}$ and $\left\{\vartheta_{-}\srr\mathfrak{Re}(\Psi_{+-})\right\}$.

The DFT spectra of elements of dictionaries $\left\{\vartheta_{+}\right\}$ and $\left\{\vartheta_{-}\right\}$ form various tiling from the pairs of quadrants $\mathbf{Q}_{0}\bigcup\mathbf{Q}_{3}$  and $\mathbf{Q}_{1}\bigcup\mathbf{Q}_{2}$ (see \eh{quadr}), respectively, by squares of different  size depending on the decomposition  level. The waveforms  shapes are close to  windowed cosines with multiple frequencies oriented in multiple directions ($2(2^{m+1}-1)$ directions at the level $m$). Combinations of waveforms from the sets $\left\{\vartheta_{+}\right\}$ and $\left\{\vartheta_{-}\right\}$ provide a variety of frames in the space of 2D signals. The  transforms are executed in a fast way using FFT.

In this paper, we design the directional qWPs starting from the discrete-time-spline  WPs. They offer more flexibility compared to the discrete-spline   WPs. In particular, the former WPs can originate from polynomial  splines of any order, while the latter WPs use only even-order discrete splines.

The paper is organized as follows:
Section \ref{sec:s2}  briefly  outlines the orthonormal  WPs originated from polynomial  splines and the corresponding transforms that serve as  a basis for the design of qWPs.
  Section \ref{sec:s3}  presents the design of qWPs
  and Section \ref{sec:s4} describes implementation of the  transforms.
   Section \ref{sec:s5} extends the design of 1D qWPs to 2D case and Section \ref{sec:s6} describes the implementation of the  transforms.
  Section \ref{sec:s7}  discusses the results and the Appendix provides proofs for two propositions.

\paragraph{Notations and abbreviations:}
$N=2^{j}$,  $\omega\srr e^{2\pi\,i/N}$ and $\Pi[N]$ is a space of real-valued  $N$-periodic  signals.
$\Pi[N,N]$ is the space of two-dimensional   $N$-periodic  arrays in both vertical and horizontal directions. The sequence  $\delta[k]\in\Pi[N]$ means the $N$-periodic  Kronecker delta.

Discrete Fourier transform  (DFT), Fast Fourier transform  (FFT), DFT of a signal  $\mathbf{x}\in\Pi[N]$ is $\hat{x}[n]=\sum_{k=0}^{N-1}\omega^{-kn}x[k]$ and $\hat{x}[n]_{m}\srr\sum_{k=0}^{2^{-m}N-1}\omega^{-2^{m}kn}x[k]$.
 $\cdot^{\ast}$ means  complex conjugate.
WPT means  wavelet packet transform  (WPT), perfect reconstruction (PR),  Hilbert transform  (HT),  $H(\mathbf{x})$ is the discrete  periodic  HT of a signal  $\mathbf{x}$.
discrete-time spline  (DTS),  DTSWP, cWP and qWP mean discrete-time-spline-based wavelet packets $\psi^{p}_{[m],l}$), complimentary wavelet packets $\varphi^{p}_{[m],l}$ and quasi-analytic  wavelet packets $\Psi^{p}_{\pm[m],l}$, respectively, in 1D case, and  wavelet packets $\psi^{p}_{[m],j,l}$, complimentary wavelet packets $\varphi^{p}_{[m],j,l}$ and quasi-analytic  wavelet packets $\Psi^{p}_{+\pm[m],l,j}$, respectively, in 2D case.
%The abbreviation {qWPA} designates the  image denoising  algorithm  based on the qWP transforms.

%Peak Signal-to-Noise ratio (PSNR) in decibels (dB) is
%\(
%  PSNR\srr10\log_{10}\left(\frac{N\,255^2}{\sum_{k=1}^N(x_{k}-\tilde
%  x_{k})^2}\right)\; dB.
%  \)
%
%SSIM means Structural Similarity Index.
%BSA stands for Bivariate Shrinkage algorithm   and
p-filter  means periodic  filter.
%SBI  stands for split Bregman iteration.
%\textbf{M1} stands for Method1, \textbf{M2} stands for Method2.
%SET-4 means the set of  the  filter banks  DAS-2,  DAS-1, TP-$\mathbb{C}$TF$_6$ and  TP-$\mathbb{C}$TF$_6^{\downarrow}$.

Notation  $l_{0}\srr0,\;l_{m}\srr2^{m}-1.$
Quadrants of the frequency  domain:
\begin{equation}\label{quadr}
\begin{array}{cc}
  \mathbf{Q}_{0}\srr[0,N/2-1]\times[0,N/2-1], & \mathbf{Q}_{1}\srr[0,N/2-1]\times[-N/2,-1], \\
 \mathbf{Q}_{2}\srr[-N/2,-1]\times[0,N/2-1], & \mathbf{Q}_{3}\srr[-N/2,-1]\times[-N/2,-1].
\end{array}
\end{equation}

\section{Outline of orthonormal WPs originated from discrete-time splines: preliminaries}\label{sec:s2}
This section provides a brief  outline of  periodic  discrete-time wavelet packets (DTSWPs) originated from polynomial  splines and corresponding transforms. For details see Chapter 4 in \cite{ANZ_book3}.
\subsection{ Periodic  discrete-time splines and first-level wavelet packets}\label{sec:ss21}

The centered  $N$-periodic   polynomial  B-spline  $B^{p}(t)$ of order $p$ is an $N$-periodization of the function
\begin{equation}
  \label{bsp_expl}
 {b}^{p}(t)=\frac{1}{(p-1)!}\sum_{k=0}^{p}(-1)^{k}\,{p\choose k}\,\left(t+\frac{p}{2}-k\right)_{+}^{p-1},\quad x_{+}\srr\max\{x,0\}.
\end{equation}
The B-spline  $B^{p}(t)$ is supported on the interval $(-p/2,p/2)$ up to periodization. It is strictly positive inside this interval and  symmetric about zero, where it has its single maximum and has $p-2$ continuous derivatives. The Fourier coefficients of the B-spline  are
\begin{equation}\label{pbspF}
c_{n}({B}^{p})=\int_{-N/2}^{N/2}B^{p}(t)\,e^{-2\pi i nt/N}=\left(\frac{\sin \pi n/N}{\pi n/N}\right)^{p}.
\end{equation}
The functions $
  {S}^{p}(t) \srr \sum_{l=0}^{N-1}q[l]\,{B}^{p}[t-l],
$
 are referred to as  the order-$p$  periodic  splines. The following two sequences (Eqs. \rf{psp_char} and \rf{psp_v}) will be repeatedly used in the further presentation:
\begin{eqnarray}\label{psp_char}
   u^{p}[n]&\srr&\sum_{k=-N/2}^{N/2-1}\omega^{-kn}\,{b}^{p}\left(k\right)=\sum_{l\in\z}\,\left(\frac{\sin \pi\,(n/N+l)}{\pi\,(n/N+l)}\right)^{p}=
   \sin^{p} \frac{\pi\,n}{N}\sum_{l\in\z}\,\frac{(-1)^{lp}}{\left(\pi\,(n/N+l)\right)^{p}},\\\label{psp_v}
v^{p}[n]&\srr&\omega^{-n/2}\sum_{k=-N/2}^{N/2-1}\omega^{-kn}\,{b}^{p}\left(k+\frac{1}{2}\right)= \sin^{p} \frac{\pi\,n}{N}\sum_{l\in\z}\,\frac{(-1)^{l(p+1)}}{\left(\pi\,(n/N+l)\right)^{p}}
%\sum_{l\in\z}(-1)^{l}\,\left(\frac{\sin \pi\,(n/N+l)}{\pi\,(n/N+l)}\right)^{p}.
\end{eqnarray}
\br\label{uv_rem}It is well known (for example, \cite{sho1}) that the $N$-periodic  sequence  $u^{p}[n]$ is strictly positive and symmetric about $N/2(\mathrm{mod}\, N) $, where it attains its single minimum. The sequence  $v^{p}[n]$ is $2N$-periodic  and $v^{p}[n+N]=-v^{p}[n]$.\er

Denote by $b_{d}^{p}(t)=b^{p}(t/2)/2$, which is  the two-times dilation of the B-spline   $b^{p}(t)$.
\bdd\label{bds_defT}The span-two discrete-time  B-spline  $\mathbf{b}_{[1]}^{p}$ of order $p$  is defined as an $N$-periodization of the sampled B-spline  $b_{d}^{p}(t)$:
$
 {b}^{p}_{[1]}[k]\srr b_{d}^{p}(k),\;k=-N/2,...,N/2-1(\mathrm{mod}\, N).
$
\edd
The discrete-time B-spline  $\mathbf{b}^{p}_{[1]}$ is an $N$-periodic  signal  from $\Pi[N]$. The DFT of the B-spline  $b^{p}_{[1]}$ is
\begin{eqnarray}\nonumber
 \hat{b}^{p}_{[1]}[n]&=&\sum_{k=-N/4}^{N/4-1}\omega^{-2kn}\,b_{d}^{p}(2k)+\omega^{-n}\,\sum_{k=-N/4}^{N/4-1}\omega^{-2kn}\,b_{[d]}^{p}(2k+1)\\\label{bts2_uv} &=&\frac{1}{2}\sum_{k=-N/4}^{N/4-1}\omega^{-2kn}\,b^{p}(k)+ \frac{\omega^{-n}}{2}\sum_{k=-N/4}^{N/4-1}\omega^{-2kn}\,b^{p}\left(k+\frac12\right) =
 \frac{u^{p}[2n]+v^{p}[2n]}{2}.
 \end{eqnarray}
The sequences  $u^{p}[n]$ and  $v^{p}[n]$  are defined in Eqs. \rf{psp_char} and \rf{psp_v}. The samples of B-splines $b^{p}(t)$ of different orders at  points $\left\{k\right\}$ $\left\{k+1/2\right\}$ can be easily computed using \eh{bsp_expl}. We used here the fact that $b^{p}(t)$ is supported on the interval
 $(-p/2,p/2)\subset (-N/4,N/4-1)\bigcup (-N/2,N/2-1)$ .

\br\label{uv2n_rem} Referring to Remark \ref{uv_rem},  we claim  that
\begin{equation}\label{b0n2}
\begin{array}{l}
  u^{p}[2n+N)]=u^{p}[2n],\;  v^{p}[2n+N)]=-v^{p}[2n], \quad u^{p}[0]=v^{p}[0]=u^{p}[N]=1,\;  \\
v^{p}[N]=-1,\quad  \hat{b}^{p}_{[1]}[0]=1, \quad\hat{b}^{p}_{[1]}[N/2]=0.
\end{array}
 \end{equation}
\er
Linear combinations of two-sample shifts of the B-splines $
    {s}_{[1]}^{p}[k]=\sum_{l=0}^{N/2-1}q[l]\,{b}_{[1]}^{p}[k-2l]$
are referred to as periodic  discrete-time splines (DTSs) of span 2. Their DFT is $\hat{{s}}_{[1]}^{p}[n]=\hat{q}[n]_{1}\,      \hat{{b}}^{p}_{[1]}[n]. $ The $N/2$-dimensional  space of the DTSs is denoted by $^{p}{{\mathcal{S}}}_{[1]}^{0}\subset \Pi[N]$.
Two-sample shifts of the discrete-time B-spline   $\mathbf{b}_{[1]}^{p}$ form a basis in the space  $^{p}{{\mathcal{S}}}_{[1]}^{0}$.  Denote by   ${} ^{p}{\mathcal{S}}_{[1]}^{1}$  the  orthogonal  complement of the subspace   ${} ^{p}{\mathcal{S}}_{[1]}^{0}$  in the signal   space   $\Pi[N]$.  Thus, $\Pi[N]={} ^{p}{\mathcal{S}}_{[1]}^{0}\bigoplus^{p}{\mathcal{S}}_{[1]}^{1}$.

Define the  DTS  ${\psi}_{[1],0}^{p}$ and the signal  ${\psi}_{[1],1}^{p}\in \Pi[N]$ by their  DFTs:
\begin{eqnarray}\label{on_dsDFT}
      \hat{{ \psi}}_{[1],0}^{p}[n]&=& \frac{\hat{b}^{p}_{[1]}[n]}{\sqrt{\Upsilon^{p}[n]}}\srr\beta[n],\quad   \hat{{ \psi}}_{[1],1}^{p}[n]= \omega^{n}\,\frac{\hat{b}^{p}_{[1]}[n+N/2]}{\sqrt{\Upsilon^{p}[n]}}\srr\alpha[n],\\\nonumber
      \Upsilon^{p}[n]&\srr &\frac{u^{p}[2n]^{2}+  v^{p}[2n]^{2}}{4}.
     \end{eqnarray}
   The real-valued signals  ${\psi}_{[1],0}^{p}[k]$ and  ${\psi}_{[1],1}^{p}[k]$ are symmetric about $k=0$ and $k=-1$, respectively.
     %Referring to Remark \ref{uv_rem},  we claim  that the sequence  $\Upsilon^{p}[n]$ is strictly positive.
   \bpp[\cite{ANZ_book3}, Chapters  3 and 4]\label{ds_proj_pro1T} Two-sample shifts of the signals   ${\psi}_{[1],\lambda}^{p}[k],\;\lambda=0,1$  form  orthonormal  bases of  the subspaces  $\,{}^{p}{{\mathcal{S}}}_{[1]}^{\lambda},\;\lambda=0,1$, respectively, such that their inner products in the space $\Pi[N]$ are $ \left\langle{\psi}_{[1],\lambda}^{p}[\cdot-2l], {\psi}_{[1],\lambda}^{p}[\cdot-2m]\right\rangle=\delta(l-m),\;\lambda=0,1.$

The orthogonal  projections of a signal   $\mathbf{x}\in \Pi[N]$  onto the subspaces   ${} ^{p}{\mathcal{S}}_{[1]}^{\lambda}$are  the signals $\mathbf{{x}}_{[1]}^{\lambda}\in\Pi[N]$, respectively, such that
\begin{eqnarray*}\nonumber
          {x}_{[1]}^{\lambda}[k]&=&\sum_{l=0}^{N/2-1}y_{[1]}^{\lambda}[l]\, {\psi}_{[1],\lambda}^{p}[k-2l] =\sum_{l=0}^{N/2-1}y_{[1]}^{\lambda}[l]\, {h}_{[1]}^{\lambda}[k-2l],\\\label{y0_delT} y_{[1]}^{\lambda}[l] &=&\left\langle \mathbf{x},\, {\psi}_{[1],\lambda}^{p}[\cdot-2l]   \right\rangle
        = \sum_{k=0}^{N-1}{h}_{[1]}^{\lambda}[k-2l] \,x[k], \quad  {h}_{[1]}^{\lambda}[k] ={\psi}_{[1],\lambda}^{p}[k],\;\lambda=0,1,\;k\in\z.
         \end{eqnarray*}
         \epp

 \br\label{filt11_remT}The sets $\left\{ y_{[1]}^{0}[l]  \right\}$ and $\left\{ y_{[1]}^{1}[l]  \right\},\;l=0,...,N/2-1,$ of the orthogonal  projection coefficients can be regarded as  results of p-filtering the signal  $ \mathbf{x}$
          by the time-reversed low- and high-pass p-filters $\mathbf{h}_{[1]}^{0}$ and $\mathbf{h}_{[1]}^{1}$, respectively, which is followed by downsampling of factor 2.
          The impulse responses of the  p-filters $\mathbf{h}_{[1]}^{j},\;\lambda=0,1,$  coincide with the signals $\psi_{[1],\lambda}^{p}[k] $, respectively. Their frequency response s  are ${h}_{[1]}^{0}=\beta[n]$, ${h}_{[1]}^{1}=\alpha[n]$. \er

          \bdd The signals $  \psi_{[1],0}^{p}$ and $  {\psi}_{[1],1}^{p}$  are referred to as the discrete-time-spline    wavelet packets (DTSWPs) of order $p$ from the first decomposition  level.\edd

\begin{SCfigure}
\centering
\caption{ Left:  DTSWPs  $ { \psi}_{[1],0}^{p}$ (red lines) and $  {\psi}_{[1],1}^{p}$  (blue lines), $p=3,8,15$.
Right: magnitude spectra of  $ { \psi}_{[1],0}^{p}$ (red lines) and $  {\psi}_{[1],1}^{p}$  (blue lines)}
\includegraphics[width=3in]{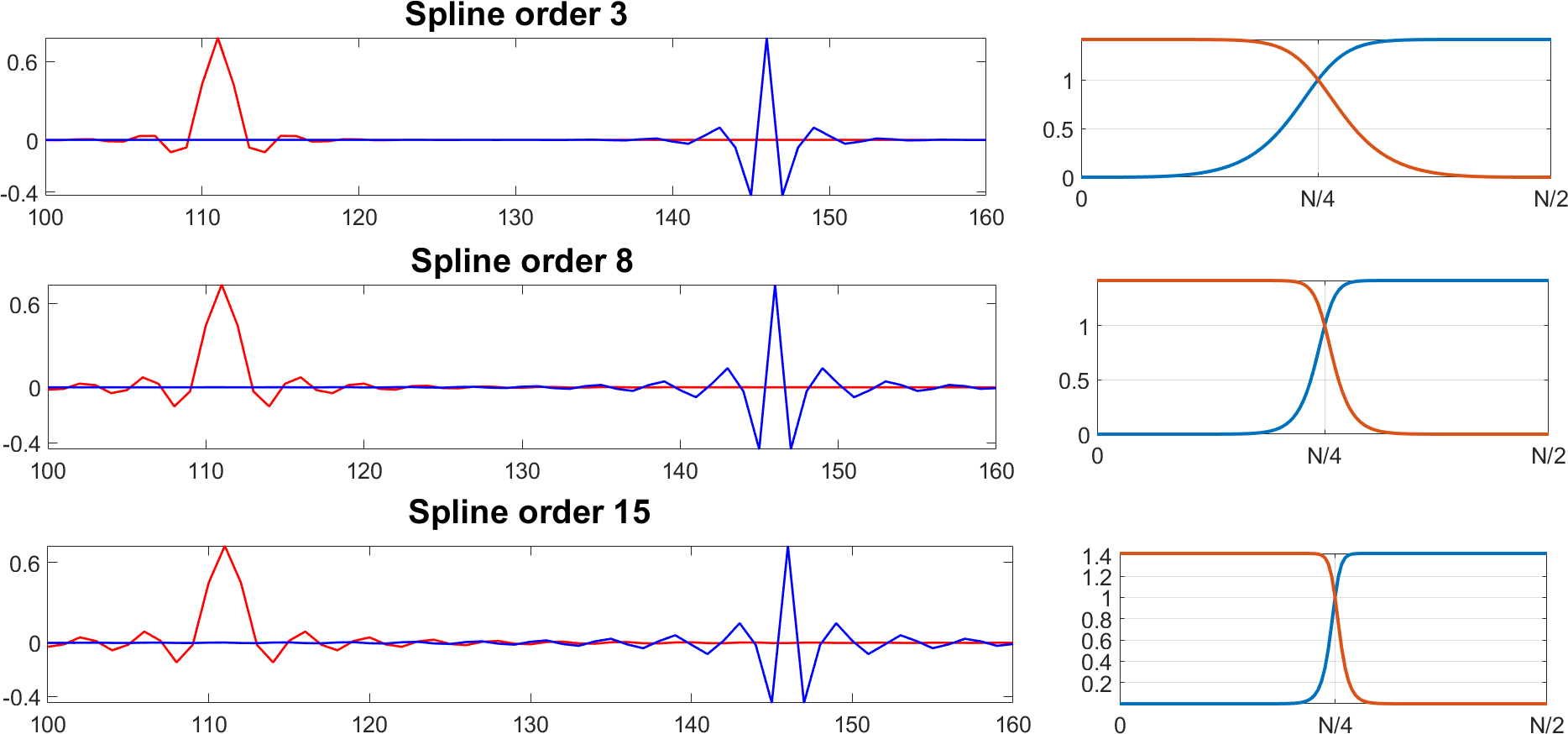}
  \label{ds_wav1LT}
\end{SCfigure}

Figure \ref{ds_wav1LT} displays the DTSWPs  $  \psi_{[1],0}^{p}$ and $  {\psi}_{[1],1}^{p}$  (which are  the  p-filters' $\mathbf{h}_{[1]}^{0}$ and $\mathbf{h}_{[1]}^{1}$  impulse responses)
and magnitudes of their DFTs (which are the p-filters' $\mathbf{h}_{[1]}^{0}$ and $\mathbf{h}_{[1]}^{1}$ magnitude responses) of different orders.
It is seen that the WPs are well localized in time domain. Their spectra are flat and their shapes tend to rectangular as their orders increase.
The one-level DTSWP transform  of a signal  $\mathbf{x}$ and its inverse are represented in a matrix form:
   \begin{equation}\label{mod_repAS1}
    \left(
     \begin{array}{c}
       \hat{y}_{[1]}^{0}[n]_{1} \\
         \hat{y}_{[1]}^{1}[n]_{1}\\
     \end{array}
   \right)=\frac{1}{2}\tilde{\mathbf{M}}[-n]\cdot \left(
     \begin{array}{l}
      \hat{x}[n] \\
       \hat{x}[\vec{n}]
     \end{array}
   \right),\quad  \left(
     \begin{array}{l}
      \hat{x}[n] \\
       \hat{x}[\vec{n}]
     \end{array}
   \right)={\mathbf{M}}[n]\cdot \left(
     \begin{array}{c}
        \hat{y}_{[1]}^{0}[n]_{1} \\
         \hat{y}_{[1]}^{1}[n]_{1}\\
     \end{array}
   \right),
    \end{equation}
    where $\vec{n}=n+{N}/{2}$ and $\tilde{\mathbf{M}}[n]$ and ${\mathbf{M}}[n]$ are the  modulation matrices of the analysis  and synthesis  p-filter banks, respectively.

The modulation matrices are:
  \begin{eqnarray}\label{sa_modma10T}
    {\mathbf{{M}}}[n]=\sqrt{2}\left(
                                 \begin{array}{cc}
                                   {\beta}[n] & {\alpha}[n]   \\
                             {\beta}\left[n+\frac{N}{2}\right]     & {\alpha}\left[n+\frac{N}{2}\right]  \\
                                 \end{array}
                               \right)= \tilde{\mathbf{{M}}}[n]^{T},
  \end{eqnarray}
  where $\beta[n]$ and $\alpha[n]$ are defined in \eh{on_dsDFT}.
  The synthesis  p-filter bank  $\mathbf{H}_{[1]}=\mathbf{h}_{[1]}^{0}\bigcup\mathbf{h}_{[1]}^{1}$  coincides with the analysis  p-filter bank  and, together, they  form a perfect  reconstruction (PR)\index{perfect  reconstruction (PR) p-filter bank} p-filter bank.

\subsection{Extension of transforms to deeper decomposition  levels}\label{sec:ss22}
\subsubsection{Second-level wavelet packet transforms (WPTs)}\label{sec:sss221}The WPT from the first to the second  decomposition  level is implemented by  application of the analysis  p-filter bank  $\tilde{\mathbf{H}}_{[2]}=\left\{\mathbf{h}^{0}_{[2]},\mathbf{h}^{1}_{[2]} \right\}$, which operates in the space $\Pi[N/2]$ to the signals
 $\mathbf{y}_{[1]}^{\lambda},\;\lambda=0,1,$.  The frequency response s of the p-filters are
 $ \hat{{h}}_{[2]}^{\lambda}[n]_{1}= \beta[2n] \quad
   \hat{{h}}_{[2]}^{1}[n]_{1} = \alpha[2n],$
where $\beta[n]$ and $\alpha[n]$ are defined in \eh{on_dsDFT}.
The modulation matrices of the     p-filter bank
${\mathbf{H}}_{[2]}$ are
 \begin{eqnarray}
 \tilde{\mathbf{M}}_{[2]}[n]=\tilde{\mathbf{M}}[2n],\quad
    {\mathbf{M}}_{[2]}[n]={\mathbf{M}}[2n]\label{sa_modma20},
  \end{eqnarray}
  where the modulation matrices $\tilde{\mathbf{M}}[n]$ and ${\mathbf{M}}[n]$ are defined in \eh{sa_modma10T}.

Define the  signals  ${\psi}_{[2],\rho}^{p}\in\Pi[N]$  by their DFT
 \begin{equation}\label{spec_psi20}
 \hat{ {\psi}}_{[2],\rho}^{p}[n]=\hat{{\psi}}_{[1],\lambda}^{p}[n]\,\hat{h}_{[2]}^{\mu}[n]_{1}=\hat{{\psi}}_{[1],\lambda}^{p}[n]\,\hat{{\psi}}_{[1],\mu}^{p}[2n]_{1},\quad \rho=\left\{
     \begin{array}{ll}
     \mu, & \hbox{if $\lambda=0$;} \\
    3-\mu, & \hbox{if $\lambda=1$.}
   \end{array}
  \right..
 \end{equation}

 \bpp[\cite{ANZ_book3}, Chapter 4]\label{psi20_pro}  The norms of the signals ${\psi}_{[2],\rho}^{p}\in \Pi[N]$ are equal to one. The 4-sample shifts $\left\{\psi_{[2],\rho}^{p}[\cdot-4l] \right\},\;l=0,...,N/4-1,$ of this signal  are mutually orthogonal  and signals with different  indices $\rho$ are orthogonal  to each other.\epp
 Thus, the signal  space $\Pi[N]$  splits into four mutually orthogonal  subspaces $\Pi[N]=\bigoplus_{\rho=0}^{3} {}^{p}\mathcal{S}_{[2]}^{\rho}$ whose orthonormal  bases are formed by 4-sample shifts  $\left\{\psi_{[2],\rho}^{p}[\cdot-4l] \right\},\;l=0,...,N/4-1,$ of the signals $\psi_{[2],\rho}^{p}$, which are referred to as the second-level DTSWPs of order $p$.

 The orthogonal   projection of a signal  $\mathbf{x}\in \Pi[N]$ onto the subspace  ${}^{p}\mathcal{S}_{[2]}^{\rho}$ is the signal
 \begin{equation*}\label{s20_rep}
 x_{[2]}^{\rho}[k]= \sum_{l=0}^{N/4-1}\left\langle \mathbf{x}, \,\psi_{[2],\rho}^{p}[\cdot-4l] \right\rangle\,{\psi}_{[2],\rho}^{p}[k-4l]=\sum_{l=0}^{N/4-1}y_{[2]}^{\rho}[l]\,{\psi}_{[2],\rho}^{p}[k-4l], ~~
 k=0, \ldots , N -1.
\end{equation*}

   Practically, derivation of the wavelet packet transform   coefficients $\mathbf{y}_{[1]}^{\lambda},\;\lambda=0,1,$ from $\mathbf{x}$ and the inverse operation are implemented using \eh{mod_repAS1}, while the transform
 $\mathbf{y}_{[1]}^{\lambda}\longleftrightarrow\mathbf{y}_{[2]}^{\rho}$ are implemented similarly using the  modulation matrices of the     p-filter bank
${\mathbf{H}}_{[2]}$ defined in \eh{sa_modma20}.
  The second-level  wavelet packets $\psi_{[2],\rho}^{p}$ are derived from the first-level  wavelet packets $\psi_{[1],\lambda}^{p}$ by filtering the latter with the p-filters  $\mathbf{h}_{[2]}^{\mu},\;\lambda,\mu=0,1$.

\begin{SCfigure}
\centering
\caption{ Left: second-level DTSWPs of different orders; left to right: $\psi_{[2],0}^{p}\to\psi_{[2],1}^{p}\to\psi_{[2],2}^{p}\to\psi_{[2],3}^{p}$. Right: Their magnitude DFT spectra}
\includegraphics[width=3in]{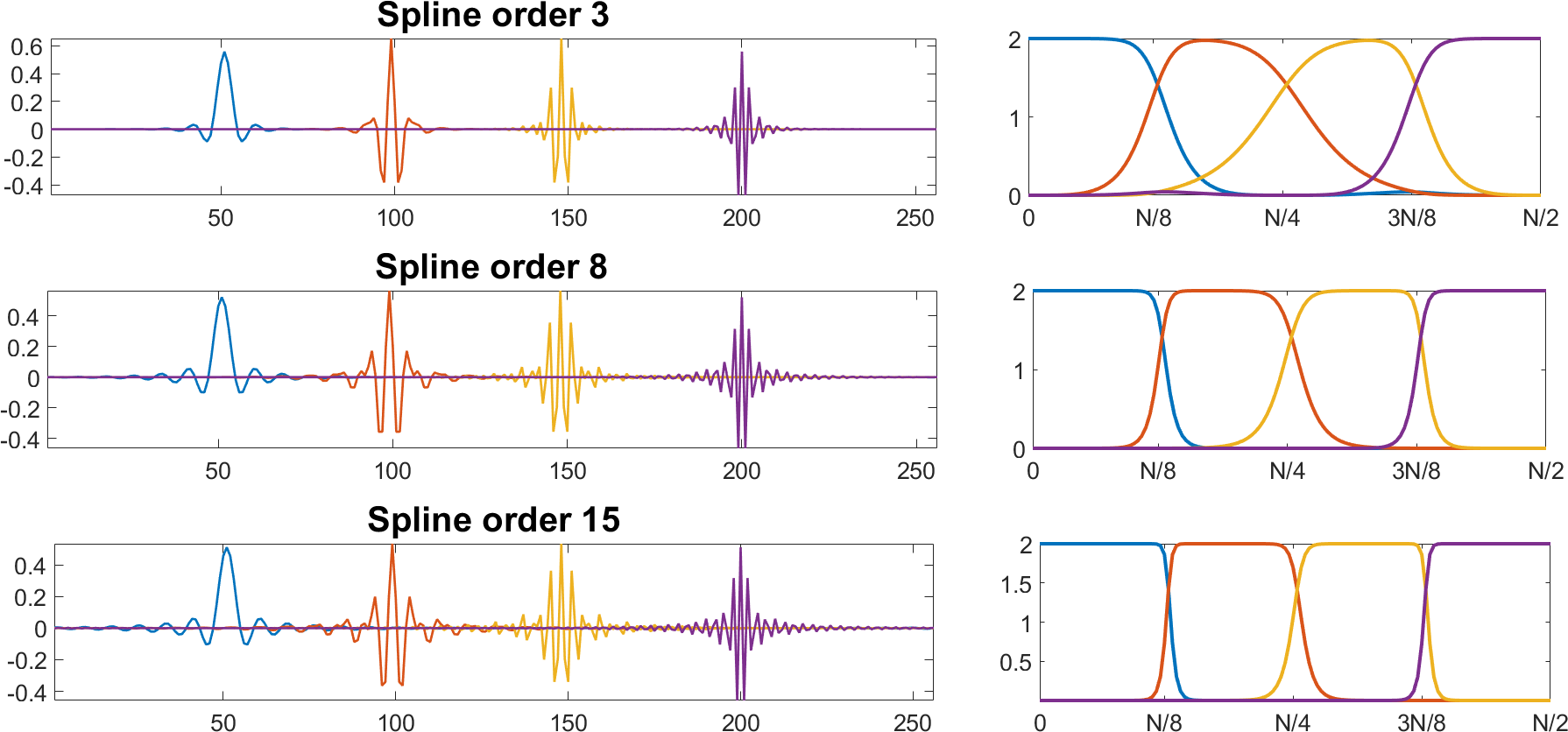}
  \label{dss_wq_s2T3}
\end{SCfigure}

 Figure \ref{dss_wq_s2T3} displays the second-level WPs originating from DTSs  of orders 3, 8 and  15   and their  DFTs. One can observe that the wavelet packets are symmetric and well localized  in time domain.
    Their spectra are flat and their shapes tend to rectangular as their orders increase. They split the frequency  domain into four quarter-bands.

\subsubsection{Transforms to deeper levels}\label{sec:sss222}
The WPTs to deeper decomposition  levels are implemented iteratively,  while the transform  coefficients $\left\{\mathbf{y}_{[m+1]}^{\rho}\right\}$ are derived by filtering the coefficients $\left\{\mathbf{y}_{[m]}^{\lambda}\right\}$ with the p-filters $\mathbf{h}^{\mu}_{[m+1]},$ where $\lambda=0,...,2^m-1,\;\mu=0,1$ and $\rho=\left\{
                                                                                                                               \begin{array}{ll}
                                                                                                                                 2\lambda +\mu, & \hbox{if $\lambda$ is even;} \\
                                                                                                                                 2\lambda +(1-\mu), & \hbox{if $\lambda$ is odd.}
                                                                                                                               \end{array}
                                                                                                                             \right.
$
The transform  coefficients are ${y}_{[m]}^{\lambda}[l]=\left\langle \mathbf{x},\psi^{p}_{[m],\lambda}[\cdot, -2^{m}l] \right\rangle$, where the signals $\psi^{p}_{[m],\lambda}$ are normalized, orthogonal  to each other in the space $\Pi[N]$, and their $2^{m}l-$sample shifts are mutually orthogonal. They are referred to as level-$m$ DTSWPs of order $p$. The set $\left\{\psi^{p}_{[m],\lambda}[\cdot, -2^{m}l] \right\},\;\lambda=0,...,2^m-1,\;l=0,...N/2^m-1,$ constitutes an orthonormal  basis of the  space $\Pi[N]$ and generates its split into $2^m$ orthogonal  subspaces. The next-level wavelet packets $\psi^{p}_{[m+1],\rho}$ are derived by filtering the wavelet packets  $\psi^{p}_{[m],\lambda}$ with the p-filters $\mathbf{h}^{\mu}_{[m+1]}$ such that
\begin{equation}\label{mlev_wq}
  {\psi}_{[m+1],\rho}^{p}[n]  =\sum_{k=0}^{N/2^{m}-1}{h}_{[m+1]}^{\mu}[k] \, {\psi}_{[m],\lambda}^{p}[n-2^{m}k].
\end{equation}
Note that the frequency response of an $m-$level p-filter  is $ \hat{h}^{\mu}_{[m]}[n]=\hat{h}^{\mu}_{[1]}[2^{m-1}n].$

The transforms are executed in the spectral domain using the  Fast Fourier transform  (FFT) by the application of critically sampled two-channel filter banks to the half-band spectral  components of a signal. For example, the Matlab execution of the 8-level 13-th-order WPT of a signal  comprising  245760 samples, takes  0.2324  seconds.

\subsection{ 2D WPTs}\label{sec:ss23}
A standard way to extend the one-dimensional (1D) WPTs to multiple dimensions is the tensor-product extension.
The 2D one-level WPT of a  signal  $\mathbf{x}=\left\{x[k,n]\right\},\;k,n=0,...,N-1,$  which  belongs to $\Pi[N,N]$, consists of the application of  1D WPT to columns of  $\mathbf{x}$, which is followed by the application of the transform  to  rows of the  coefficient  array. As a result of the 2D WPT of signals from $\Pi[N,N]$, the space becomes split
 into four mutually orthogonal  subspaces
$ \Pi[N,N]=\bigoplus_{j,l=0}^{1}\,^{p} {\mathcal{S}}^{j,l}_{[1]}.$

The 2D wavelet packets are $\psi_{[1],j,l}^{p}[n,m]\srr \psi_{[1],j }^{p}[n]\, \psi_{[1],l}^{p}[m],   \quad j,l=0,1.$ They are normalized and  orthogonal  to each other in the space $\Pi[N,N]$. It means that \\ $\sum_{n,m=0}^{N-1}\psi_{[1],j1 ,l1}^{p}[n,m]\,\psi_{[1],j2 ,l2}^{p}[n,m]=\delta[j1-j2]\,\delta[l1-l2]$. Their two-sample shifts in both directions are mutually orthogonal.
The subspace ${}^{p} {\mathcal{S}}^{j,l}_{[1]}$ is a linear hull of two-sample shifts of the 2D wavelet packets
$\left\{\psi_{[1],j,l}^{p}[k-2p,n-2t]\right\} ,\;p,t,=0,...,N/2-1,$ that form an orthonormal  basis of ${}^{p} {\mathcal{S}}^{j,l}_{[1]}$. The orthogonal  projection of the  signal  $\mathbf{x}\in\Pi[N,N]$ onto the subspace ${}^{p} {\mathcal{S}}^{j,l}_{[1]}$ is the  signal  $\mathbf{x}_{[1]}^{j,l}\in\Pi[N,N]$ such that
%\begin{equation*}\label{op_2d}
$ {x}_{[1]}^{j,l}[k,n]=\sum_{p,t=0}^{N/2-1} y_{[1]}^{j,l}[p,t] \,\psi_{[1],j,l}^{p}[k-2p,n-2t], \quad j,l =0,1.$  The transform  coefficients are $$y_{[1]}^{j,l}[p,t] =\left\langle \mathbf{x},\psi_{[1],j,l}^{p}[\cdot-2p,\cdot-2t] \right\rangle=\sum_{n,m=0}^{N-1} \psi_{[1],j,l}^{p}[n-2p,m-2t]\: x[n,m].$$

By the application of the above transforms iteratively to blocks of the transform  coefficients down to $m$-th level, we get that the space $ \Pi[N,N]$ is decomposed into $4^{m}$ mutually orthogonal  subspaces
$ \Pi[N,N]=\bigoplus_{j,l=0}^{2^{m}-1}\,^{p} {\mathcal{S}}^{j,l}_{[m]}.$  The orthogonal  projection of the  signal  $\mathbf{x}\in\Pi[N,N]$ onto the subspace ${}^{p} {\mathcal{S}}^{j,l}_{[m]}$ is the  signal  $\mathbf{x}_{[m]}^{j,l}\in\Pi[N,N]$ such that
\begin{eqnarray*}\label{op_2d}
 {x}_{[m]}^{j,l}[k,l]&=&\sum_{p,t=0}^{N/2^{m}-1} y_{[m]}^{j,l}[p,t] \,\psi_{[m],j ,l}^{p}[k-2^{m}p,l-2^{m}t], \quad j,l =0,...,2^{m}-1,\\
\psi_{[m],j ,l}^{p}[k,n]&=&\psi_{[m],\lambda}^{p}[k]\,\psi_{[m],l}^{p}[n],\quad  y_{[m]}^{j,l}[p,t]
=\left\langle \mathbf{x},\psi_{[m],j ,l}^{p}[\cdot-2^{m}p,\cdot-2^{m}t] \right\rangle.
\end{eqnarray*}

The 2D tensor-product  wavelet packets $\psi_{[m],j ,l}^{p}$ are well localized in the spatial domain, their 2D DFT spectra  provide a refined split of the frequency  domain of signals from $ \Pi[N,N].$\footnote{Especially it is true for WPs derived from  higher-order DTSs.} The drawback is that the WPs are  oriented  in ether horizontal or vertical directions or are not oriented  at all.

\subsection{Local discrete  vanishing moments}\label{sec:ss24}
One of fundamental features of wavelets and wavelet packets is their vanishing moment  property. In a conventional setting it means the annihilation of polynomials of a certain degree by a continuous wavelet  or wavelet packet $\psi(t)$. To be specific, if for any polynomial  $P_{m-1}(t)$ of degree $m-1$ the relation $\int\psi(t)P_{m-1}(t)\,dt=0$ holds, then it is said that $\psi(t)$ has $m$ vanishing moments.

We modify   the  vanishing moment  property for the discrete  periodic  setting.
 \bpp[\cite{ANZ_book1}, Chapter 15]\label{dvm_proP}Assume that the
frequency response of the high(band)-pass p-filter  ${\mathbf{g}}$ can be represented
as \(
  \hat{{g}}[n]=\sin\left(\frac{\pi n}{N}\right)^{m}\,\xi[n],
\) where $m$ is some natural number, and $\xi[n]$ is  an $N-$periodic  sequence. Assume that $\mathbf{p} $ is a signal  from $\Pi[N]$, and it
coincides with a sampled polynomial  $\mathbf{P}_{m-1}$ of degree $m-1$ at
some interval $p[k]={P}_{m-1}(k)$ as $k=k_{0},...,k_{m}$, where
$m<k_{m}-k_{0}<N$. Then, \(
  \sum_{l=0}^{N-1}g[k-l]\,p[l]=0, \mbox{   as } k=k_{0},...,k_{m}-m-1.
\) \epp

 \bdd\label{dvm__defP} If a high(band)-pass p-filter 
 ${\mathbf{g}}$  satisfies  the conditions of  Proposition \ref{dvm_proP},
we say that the p-filter  ${\mathbf{g}}$ locally eliminates sampled polynomials
of degree $m-1$. If a wavelet packet is $\psi^{p}_{[l],j}[k]\srr g[k],\;k\in\z,$ we say
that the wavelet packet $\psi^{p}_{[l],j}$ has $m$ local discrete  vanishing moments (LDVMs).\edd

\bpp\label{LDVM_pro}Assume that $\psi^{p}_{[l],j},\;j=1,...,2^{l}-1$ is a DTSWP from the decomposition   level $l$, which is  derived from the spline of order $p$. If $p$ is equal to either $2r-1$ or $2r$, then the wavelet packet $\psi^{p}_{[l],j}$ has $2r$ LDVMs.
\epp
\proof In Appendix.

\section{(Quasi-)analytic and complementary WPs}\label{sec:s3}
In this section, we define analytic  and the so-called quasi-analytic   WPs related to the DTSWPs discussed in Section \ref{sec:s2} and introduce an orthonormal  set of waveforms which are complementary to the above WPs.
\subsection{ Analytic periodic  signals}\label{sec:ss31}

A signal  $\mathbf{x}\in\Pi[N]$ is represented by its  inverse DFT which can be written as follows:
\begin{eqnarray*}\label{fsx+}
  x[k]&=&\frac{\hat{x}[0]+(-1)^{k}\hat{x}[N/2]}{N}+\frac{2}{N}\sum_{n=1}^{N/2-1}\frac{\hat{x}[n]\,\omega^{kn}+(\hat{x}[n]\,\omega^{kn})^{\ast}}{2}.
\end{eqnarray*}
Define the real-valued  signal  $\mathbf{h}\in\Pi[N]$ and two complex-valued signals $\mathbf{{x}}_{+}$ and  $\mathbf{{x}}_{-}$ such that
\begin{equation}
\label{yy}
\begin{array}{lll}
 h[k]&\srr&\frac{2}{N}\sum_{n=1}^{N/2-1}\frac{\hat{x}[n]\,\omega^{kn}-\hat{x}[n]^{\ast}\,\omega^{-kn}}{2i},\\
 {x}_{\pm}[k]&\srr&x[k]\pm ih[k]=\frac{\hat{x}[0]+(-1)^{k}\hat{x}[N/2]}{N}\\&+&\frac{2}{N}\sum_{n=1}^{N/2-1}
 \left\{
   \begin{array}{ll}
     \hat{x}[ n]\,\omega^{ kn}, & \hbox{for $\bar{x}_{+}$;} \\
      \hat{x}[- n]\,\omega^{- kn}=\hat{x}[ N-n]\,\omega^{- k(N-n)}, & \hbox{for $\bar{x}_{-}$.}
   \end{array}
 \right.
\end{array}
\end{equation}

The spectrum of $\mathbf{{x}}_{+}$ comprises only non-negative frequencies and vice versa for $\mathbf{{x}}_{-}$. We have
$\mathbf{x}=\mathfrak{Re}(\mathbf{{x}}_{\pm})$ and $\mathfrak{Im}(\mathbf{{x}}\pm)=\pm\mathbf{h}$. The  signals  $\mathbf{{x}}_{\pm}$ are referred to as  periodic  analytic signals.

Thus, the signal  $\mathbf{h}$    can be regarded as a discrete  periodic  version of the  Hilbert transform  (HT) of a discrete-time  periodic  signal   $\mathbf{x}$, that is  $\mathbf{h}=H(\mathbf{x})$ (see \cite{opp}, for example).
\bpp\label{pro:ysym}
\par\noindent
\begin{enumerate}
  %\item The HT  $\mathbf{h}=H(\mathbf{x})$  is invariant with respect to circular shift in $\Pi[N]$. That means that $\mathbf{\tilde{h}}=\mathbf{h}[\cdot +m]$ is the HT of  $\mathbf{\tilde{x}}=\mathbf{x}[\cdot +m]$.
      \item  If   the signal  $\mathbf{x}\in\Pi[N]$  is symmetric about a grid point $k=K$ than   $\mathbf{h}=H(\mathbf{x})$ is antisymmetric about \emph{K}  and $h[K]=0$.
 \item Assume that a   signal  $\mathbf{x}\in\Pi[N]$ and $\hat{x}[0]=\hat{x}[N/2]=0$. Then,
\begin{enumerate}
  \item The norm of  its HT  is $\|H(\mathbf{x}) \|=\|\mathbf{x} \|$.
  \item The magnitude spectra of the signals $\mathbf{x}$ and $\mathbf{h}=H(\mathbf{x})$ coincide.
  \end{enumerate}
\end{enumerate}\epp
\proof straightforward.

\subsection{Analytic  WPs}\label{sec:ss32}
Denote $l_{0}\srr0,\;l_{m}\srr2^{m}-1.$

The analytic  DTSWPs and their DFT spectra are  derived from the corresponding DTSWPs $\left\{\psi^{p}_{[m],l}\right\},\;m=1,...,M,\;l=0,...,2^{m}-1,$ in line with the scheme in Section \ref{sec:ss31}. Recall that for all $l\neq l_{0}$, the DFT $\hat{\psi}^{p}_{[m],l}[0]=0$ and  for all $l\neq l_{m}$, the DFT $\hat{\psi}^{p}_{[m],l}[N/2]=0$.

Denote by $\theta^{p}_{[m],l}=H(\psi^{p}_{[m],l})$ the HT  of the wavelet packet $\psi^{p}_{[m],l}$, such that the DFT is\\
%\begin{equation*}\label{th_df}
$  \hat{\theta}^{p}_{[m],l}[n]=\left\{
               \begin{array}{ll}
                  -i\,\hat{\psi}^{p}_{[m],l}[n], & \hbox{if $0<n<N/2$;} \\
                i\, \hat{\psi}^{p}_{[m],l}[n], & \hbox{if $-N/2<n<0$;} \\
                 0, & \hbox{if $n=0$, or  $n=N/2$ .}
              \end{array}
             \right.$
%\end{equation*}
\par\noindent
 Then, the corresponding analytic  DTSWPs are
\(%begin{equation*}\label{awq}
{\psi}^{p}_{\pm[m],l}=\psi^{p}_{[m],l}   \pm i\theta^{p}_{[m],l}.
\)%end{equation*}
\paragraph{Properties of the analytic  WPs}
\begin{enumerate}
  \item The DFT spectra of the analytic  WPs ${\psi}^{p}_{+[m],l}$ and ${\psi}^{p}_{-[m],l}$ are located within the bands $[0,N/2]$ and $[N/2,N]\Longleftrightarrow[-N/2,0]$, respectively.
  \item The real component ${\psi}^{p}_{[m],l}$ is the same for both WPs ${\psi}^{p}_{+[m],l}$ and ${\psi}^{p}_{-[m],l}$. {It} is a symmetric oscillating waveform.
  \item\label{prop3}
The HT WPs $\theta^{p}_{[m],l}=H({\psi}^{p}_{[m],l})$ are antisymmetric oscillating waveforms.
\item For all $l\neq l_{0}, \,l_{m}$, the  norms  $\left\| \theta^{p}_{[m],l}\right\|=1$.  Their magnitude spectra $\left|\hat{\theta}^{p}_{[m],l}[n]\right|$  coincide with the magnitude spectra of the respective WPs  $\psi^{p}_{[m],l}$.
    \item\label{prop5}  When $l= l_{0}$ or  $l=l_{m}$,  the magnitude spectra of $\theta^{p}_{[m],l}$ coincide with that of ${\psi}^{p}_{[m],l}$ everywhere except for the points $n=0$ or $N/2,$ respectively,  and the waveforms'  norms are no longer equal to 1.
\end{enumerate}
Properties in items \ref{prop3}--\ref{prop5} follow directly from Proposition \ref{pro:ysym}.

\bpp\label{pro:teta_oo} { For all $l\neq l_{0}, \,l_{m}$,
the shifts of  the HT WPs $\left\{\theta^{p}_{[m],l}[\cdot-2^{m}l]\right\}$ are orthogonal  to each other in the space $\Pi[N]$. The orthogonality does not take place for  for $\theta^{p}_{[m],0}$ and $\theta^{p}_{[m],2^{m}-1}$.}
\epp
\proof Assume that  $l\neq l_{0}, \,l_{m}$. The inner product is
\begin{eqnarray*}
   && \left\langle \theta^{p}_{[m],l},\theta^{p}_{[m],l}[\cdot-2^{m}l]  \right\rangle =\frac{1}{N}\sum_{n=-N/2}^{N/2-1}\omega^{2^{m}ln}\left|\hat{\theta}^{p}_{[m],l}[n] \right|^{2}\\
&& =\frac{1}{N}\sum_{n=-N/2}^{N/2-1}\omega^{2^{m}ln}\left|\hat{\psi}^{p}_{[m],l}[n] \right|^{2}=\left\langle \psi^{p}_{[m],l},\psi^{p}_{[m],l}[\cdot-2^{m}l]  \right\rangle=0.
\end{eqnarray*}
\eop

\subsection{Complementary set of wavelet packets and quasi-analytic  WPs}\label{sec:ss33}
\subsubsection{Complementary orthonormal   WPs}\label{sec:sss331}
The values $\hat{\theta}^{p}_{[m],j}[0]$ and $\hat{\theta}^{p}_{[m],j}[N/2]$ are missing in the DFT spectra of the HT waveforms $\theta^{p}_{[m],0}$ and $\theta^{p}_{[m],2^{m}-1}$, which the set  $\left\{\theta^{p}_{[m],j}\right\}$  from forming orthonormal  bases in the corresponding subspaces.

This keeping in mind, we define  a set $\left\{\varphi^{p}_{[m],l}\right\},\;m=1,...,M, \;l=0,...,2^{m}-1,$  of signals from the space $\Pi[N]$ via their DFTs:
\begin{equation}\label{phi_df}
  \hat{\varphi}^{p}_{[m],l}[n]= \hat{\theta}^{p}_{[m],l}[n]+
                 \hat{\psi}^{p}_{[m],l}[0]+\hat{\psi}^{p}_{[m],l}[N/2].
\end{equation}
For all $l\neq l_{0}, l_{m},$ the signals $\varphi^{p}_{[m],l}$ coincide with $\theta^{p}_{[m],l}=H(\psi^{p}_{[m],l})$.
\bpp\label{pro:phi_oo}
\par\noindent
\begin{description}
  \item[-] The magnitude spectra  $\left|\hat{\varphi}^{p}_{[m],l}[n]\right|$  coincide with the magnitude spectra of the respective WPs  $\psi^{p}_{[m],l}$.
  \item[-] For any   $m=1,...,M,$ and $l=1,...,2^{m}-2,$ the signals \ $\varphi^{p}_{[m],l}$ are antisymmetric oscillating waveforms.   For $l= l_{0}, \,l_{m}$, the shapes of the signals are near antisymmetric.
      \item[-] The orthonormality properties that are similar to the properties of WPs  $\psi^{p}_{[m],l}$ hold for the signals  $\varphi^{p}_{[m],l}$ such that
      $        \left\langle\varphi^{p}_{[m],l}[\cdot -p\,2^{m}],\varphi^{p}_{[m],\lambda}[\cdot -s\,2^{m}] \right\rangle= \delta[\lambda,l]\,\delta[p,s].
     $% \end{eqnarray*}
\end{description}
\epp
The proof of Proposition \ref{pro:phi_oo} is similar to the proof of Proposition \ref{pro:teta_oo}. 

\begin{SCfigure}%[H]
\caption{Top: signals ${\psi}^{9}_{[3],l},\;l=0,...,7$. Center: signals ${\varphi}^{9}_{[3],l},\;l=0,...,7$.  Bottom:  their magnitude DFT spectra, respectively}
\centering
\includegraphics[width=4.5in]{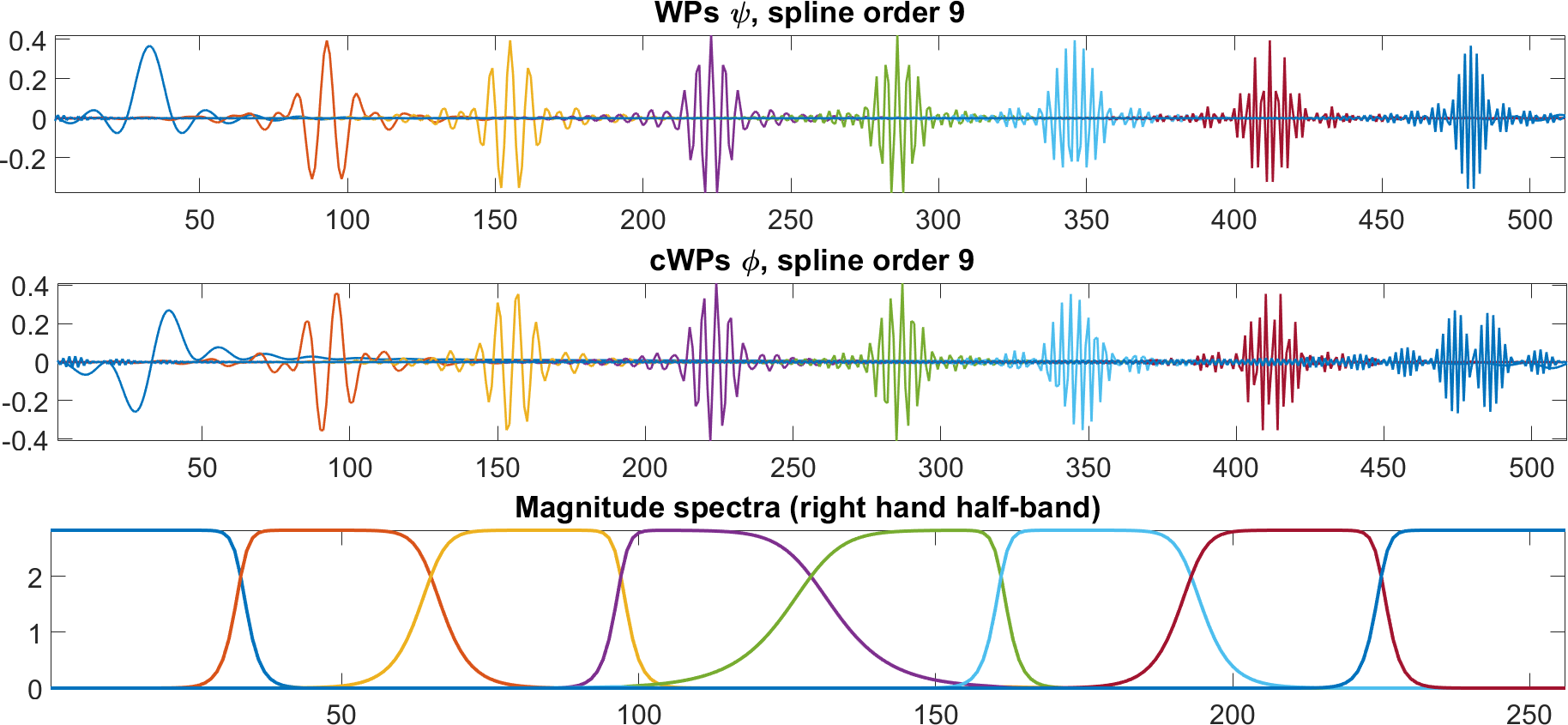}
    \label{psi_phiT}
\end{SCfigure}

  Figure \ref{psi_phiT} displays  the signals ${\psi}^{9}_{[3],l}$ and ${\varphi}^{9}_{[3],l},\;l=0,...,7$, from the third decomposition  level and their magnitude spectra.  Addition of $\hat{\psi}^{p}_{[3],l}[0]$ and $\hat{\psi}^{p}_{[3],l}[N/2]$ to the  spectra of ${\varphi}^{9}_{[3],l},\;l=0,7$ results in an antisymmetry distortion.
%\paragraph{Properties}

We call the signals     $\left\{\varphi^{p}_{[m],l}\right\},\;m=1,...,M, \;l=0,...,2^{m}-1$, the \emph{complementary  wavelet packets} (cWPs). Similarly to the DTSWPs $\left\{\psi^{p}_{[m],l}\right\},$ differentent combinations of the cWPs can provide differentent  orthonormal  bases for the space $\Pi[N]$.  These  can be, for example, the wavelet  bases
 or a   type of Best Basis \cite{coiw1,sai}.
\subsubsection{Quasi-analytic   WPs}\label{sec:sss332}
The sets of complex-valued WPs, which we refer to as the quasi-analytic   wavelet packets (qWP),  are defined by 
  $\Psi^{p}_{\pm[m],l}=\psi^{p}_{[m],l}   \pm i\varphi^{p}_{[m],l}, \quad m=1,...,M,\;l=0,...,2^{m}-1,$
where $\varphi^{p}_{[m],l}$ are the cWPs from \eh{phi_df}. The qWPs $\Psi^{p}_{\pm[m],l}$ differ from the analytic  WPs ${\psi}^{p}_{\pm[m],l}$ by the addition of the two values $\pm i\,\hat{\psi}^{p}_{[m],l}[0]$ and $\pm i\,\hat{\psi}^{p}_{[m],l}[N/2]$ into their DFT spectra, respectively. For a given decomposition  level $m$, these values are zero for all $l$ except for $l_{0}=0$ and $l_{m}=2^{m}-1$. It means that for all $l$ except for $l_{0}$ and $l_{m}$, the qWPs $\Psi^{p}_{\pm[m],l}$ are analytic.
The DFTs of qWPs are
\begin{eqnarray}\label{qa_df}
  \hat{\Psi}^{p}_{+[m],l}[n]=\left\{
               \begin{array}{ll}
  (1+i)\hat{\psi}^{p}_{[m],l}[n], & \hbox{if $n=0, N/2$;} \\
                  2\hat{\psi}^{p}_{[m],l}[n], & \hbox{if $0< n< N/2$;} \\
               0 & \hbox{if $ N/2<n<N$,}
               \end{array}
             \right.\quad
 \hat{\Psi}^{p}_{-[m],l}[n]=\left\{
               \begin{array}{ll}
(1-i)\hat{\psi}^{p}_{[m],l}[n], & \hbox{if $n=0, N/2$;} \\
                  0 & \hbox{if $0< n< N/2$;} \\
               2\hat{\psi}^{p}_{[m],l}[n], & \hbox{if $ N/2< n< N$.}
               \end{array}
             \right.
\end{eqnarray}
\subsubsection{Design of cWPs and qWPs}\label{sec:sss333}
The DFTs of the first-level DTSWPs are
%begin{equation*}\label{df_wq1}
  $\hat{\psi}^{p}_{[1],0}[n]=\frac{\hat{b}^{p}_{[1]}[n]}{\sqrt{\Upsilon^{p}[n]}}=\beta[n],\quad  \hat{\psi}^{p}_{[1],1}[n]=\omega^{n}\,\beta[n+N/2]=\alpha[n],$
%\end{equation*}
where the sequence  $ \Upsilon^{p}[n]$ is defined in \eh{on_dsDFT}. \ehh{b0n2} implies that $ \hat{\psi}^{p}_{[1],0}[0]=\sqrt{2}$ and $ \hat{\psi}^{p}_{[1],1}[N/2]=-\sqrt{2}$.

Consequently, the DFTs of  the  first-level cWPs are
\begin{eqnarray}\label{df_cwq1}
  \hat{\varphi}^{p}_{[1],0}[n]=\left\{
               \begin{array}{ll}
                  -i\,\beta[n], & \hbox{if $0<n<N/2$;} \\
                i\, \beta[n], & \hbox{if $N/2<n<N$;} \\
                 \sqrt{2}, & \hbox{if $n=0$; }
             \\
                 0, & \hbox{ if $n=N/2,$}
               \end{array}
\right.\quad \hat{\varphi}^{p}_{[1],1}[n]=\left\{
               \begin{array}{ll}
                  -i\,\alpha[n], & \hbox{if $0<n<N/2$;} \\
                i\, \alpha[n], & \hbox{if $N/2<n<N$;} \\
                 0, & \hbox{if $n=0$; }
             \\
                 -\sqrt{2}, & \hbox{ if $n=N/2.$}
               \end{array}
\right.
\end{eqnarray}
\bpp\label{pro:cwq_des} Assume that for a DTSWP $\psi_{[m+1],\rho}^{p}$ the relation in \eh{mlev_wq} holds. Then, for the cWP $\varphi_{[m+1],\rho}^{p}$ we have
\begin{eqnarray*}\label{mlev_cwq}
  {\varphi}_{[m+1],\rho}^{p}[n]  &=&\sum_{k=0}^{N/2^{m}-1}{h}_{[m+1]}^{\mu}[k] \, {\varphi}_{[m],\lambda}^{p}[n-2^{m}k]\Longleftrightarrow\hat{\varphi}_{[m+1],\rho}^{p}[\nu]=
 \hat{h}_{[1]}^{\mu}[2^{m}\nu]_{m}\,\hat{\varphi}_{[m],\lambda}^{p}[\nu],\\\nonumber
\hat{h}_{[1]}^{0}[\nu] &=&\hat{\psi}^{p}_{[1],0}[\nu]=\beta[\nu],\quad  \hat{h}_{[1]}^{1}[\nu] = \hat{\psi}^{p}_{[1],1}[\nu]=\alpha[\nu].
\end{eqnarray*}
\epp
\proof 
Due to \eh{spec_psi20}, the  DFT of the second-level DTSWPs are
\begin{eqnarray}\nonumber
 \hat{ {\psi}}_{[2],\rho}^{p}[n]&=&\hat{{\psi}}_{[1],\lambda}^{p}[n]\,\hat{h}_{[2]}^{\mu}[n]_{1},\quad \lambda,\mu=0,1,\;\rho=2\lambda+\left\{
                                                   \begin{array}{ll}
                                                     \mu, & \hbox{if $\lambda=0$;} \\
                                                     1-\mu, & \hbox{if $\lambda=1$.}
                                                   \end{array}
                                                 \right.,\\\label{df_wq2}
\hat{h}_{[2]}^{0}[n]_{1}&=&   \beta[2n],\quad \hat{h}_{[2]}^{1}[n]_{1}=\alpha[2n].
 \end{eqnarray}
For example, assume that $\lambda=\mu=0$. Then we have
\(\hat{ \psi}_{[2],0}^{p}[n]=\hat{{\psi}}_{[1],0}^{p}[n]\,\hat{h}_{[2]}^{0}[n]_{1}=\beta[n]\,\beta[2n].\)
Keeping in mind that the sequence  $\beta[2n]$  is $N/2-$periodic,  we have that the DFT of the corresponding cWP is
\[\hat{\varphi}^{p}_{[2],0}[n]=\beta[0]^2+\widehat{H({ \psi}_{[2],0}^{p})}[n]=\beta[2n]\,\left\{
               \begin{array}{ll}
                  -i\,\beta[n], & \hbox{if $0<n<N/2$;} \\
                i\, \beta[n], & \hbox{if $N/2<n<N$;} \\
                 \sqrt{2}, & \hbox{if $n=0$; }
             \\
                 0, & \hbox{ if $n=N/2,$}
               \end{array}
\right.=\hat{{\varphi}}_{[1],0}^{p}[n]\,\hat{h}_{[2]}^{0}[n]_{1}=\hat{{\varphi}}_{[1],0}^{p}[n]\,\hat{h}_{[1]}^{0}[2n]_{1}.
\]
A similar reasoning is applicable to  all the second-level cWPs and to the  cWPs from further decomposition  levels.\eop
\bcc\label{cor:cwq_des}Assume that for a DTSWP $\psi_{[m+1],\rho}^{p}$ the relation in \eh{mlev_wq} holds. Then, for the qWP $\Psi_{\pm[m+1],\rho}^{p}$ we have
\begin{eqnarray}\label{mlev_qwq}
  \Psi_{\pm[m+1],\rho}^{p}[n]  =\sum_{k=0}^{N/2^{m}-1}{h}_{[m+1]}^{\mu}[k] \, \Psi_{\pm[m],\lambda}^{p}[n-2^{m}k]\Longleftrightarrow\hat{\Psi}_{\pm[m+1],\rho}^{p}[\nu]=
 \hat{h}_{[1]}^{\mu}[2^{m}\nu]_{m}\,\hat{\Psi}_{\pm[m],\lambda}^{p}[\nu].
\end{eqnarray}
\ecc
\br\label{one_fil_rem}We emphasize that in order to derive  the $m+1$-level cWPs and qWPs from the $m$-level ones, the same p-filters are used that are used for deriving the corresponding $m+1$-level DTSWPs from the $m$-level ones.\er

\section{Implementation of cWP and qWP transforms}\label{sec:s4}
Implementation of transforms with DTSWPs ${\psi}_{[m],\lambda}^{p}$ was discussed in Section \ref{sec:s2}. In this section, we extend the transform  scheme to the transforms with cWPs ${\varphi}_{[m],\lambda}^{p}$ and qWPs $\Psi_{[m],\lambda}^{p}$.
\subsection{One-level transforms}\label{sec:ss41}
Denote by ${} ^{p}{\mathcal{C}}_{[1]}^{0}$ the subspace of the signal  space $\Pi[N]$, which is the linear hull of the set $\mathbf{W}_{[1]}^{0}=\left\{{\varphi}_{[1],0}^{p}[\cdot-2k]\right\},\;k=0,...,N/2-1$. The signals from the set $\mathbf{W}_{[1]}^{0}$ form an orthonormal  basis of the subspace  ${} ^{p}{\mathcal{C}}_{[1]}^{0}$. Denote by ${} ^{p}{\mathcal{C}}_{[1]}^{1}$ the orthogonal  complement of the  subspace ${} ^{p}{\mathcal{C}}_{[1]}^{0}$ in the  space $\Pi[N]$. The signals from   the set $\mathbf{W}_{[1]}^{1}=\left\{{\varphi}_{[1],1}^{p}[\cdot-2k]\right\},\;k=0,...,N/2-1$ form an orthonormal  basis of the subspace  ${} ^{p}{\mathcal{C}}_{[1]}^{1}$.
\bpp\label{pro:phi1}
The orthogonal  projections of a signal   $\mathbf{x}\in \Pi[N]$  onto the spaces   ${} ^{p}{\mathcal{C}}_{[1]}^{\mu},\;\mu=0,1$ are the signals $\mathbf{x}_{[1]}^{\mu}\in\Pi[N]$ such that
\begin{eqnarray*}\label{c1_del1}
          x_{[1]}^{\lambda}[k]&=&\sum_{l=0}^{N/2-1}c_{[1]}^{\lambda}[l]\, \varphi_{[1],\lambda}^{p}[k-2l] ,\quad
          %\\label{c1_del}
           c_{[1]}^{\lambda}[l]  =\left\langle \mathbf{x},\, \varphi_{[1],\lambda}^{p}[\cdot-2l]   \right\rangle
       = \sum_{k=0}^{N-1}g_{[1]}^{\lambda}[k-2l] \,x[k], \\\nonumber  g_{[1]}^{\lambda}[k] &=&\varphi_{[1],\lambda}^{p}[k],\quad
         \hat{{g}}_{[1]}^{\lambda}[n] = \hat{ \varphi}_{[1],\lambda}^{p}[n],\quad \lambda=0,1.
         \end{eqnarray*}
The DFTs $ \hat{ \varphi}_{[1],\lambda}^{p}[n]$ of the first-level cWPs are given in \eh{df_cwq1}.
\epp
The  transforms $\mathbf{x}\rightarrow \mathbf{c}_{[1]}^{0}\bigcup\mathbf{c}_{[1]}^{1}$ and back are implemented using the analysis  $ \tilde{\mathbf{M}}^{c}[n] $ and the synthesis
$ \mathbf{M}^{c}[n] $ modulation matrices:

 \begin{equation}
 \label{aa_modma11}
   \begin{array}{lll}
 \tilde{\mathbf{M}}^{c}[n]&\srr& \left(
                         \begin{array}{cc}
                            \hat{g}_{[1]}^{0}[n] &   \hat{g}_{[1]}^{0}\left[n+\frac{N}{2}\right]\\
                            \ \hat{g}^{1}_{[1]}[n] &   \hat{g}^{1}_{[1]}\left[n+\frac{N}{2}\right]\\
                         \end{array}
                       \right)=\left(\begin{array}{cc}
                                   \check{\beta}[n] & -\check{\beta}\left[n+\frac{N}{2}\right]  \\
                                  \check{\alpha}[n] &- \check{\alpha}\left[n+\frac{N}{2}\right]  \\
                                 \end{array}
                               \right),\\
    {\mathbf{M}}^{c}[n]&\srr& \left(\begin{array}{cc}
                                   \check{\beta}[n] &            \check{\alpha}[n] \\
                 -\check{\beta}\left[n+\frac{N}{2}\right]     &- \check{\alpha}\left[n+\frac{N}{2}\right]  \\
                                 \end{array}
                               \right),\quad n-0,...N/2,
\\  % \label{alphabetC}
    \check{\beta}[n]      &=&      \left\{
                                                                                                   \begin{array}{ll}
                                                                                                     {\beta}[0], & \hbox{if $n=0$;} \\
                                                                                                     -i{\beta}[n], & \hbox{otherwise,}
                                                                                                   \end{array}
                                                                                                 \right.
     \quad
       \check{\alpha}[n]= \left\{
                                                                                                   \begin{array}{ll}
                                                                                                     {\alpha}[N/2], & \hbox{if $n=N/2$;} \\
                                                                                                     -i\alpha[n], & \hbox{otherwise.}
                                                                                                   \end{array}
                                                                                                 \right.
  \end{array}
  \end{equation}
The sequences $\beta[n]$ and $\alpha[n]$ are given in \eh{on_dsDFT}.

Similarly to \eh{mod_repAS1}, the one-level cWP transform  of a signal  $\mathbf{x}$ and its inverse are:
   \begin{equation*}\label{mod_repAC1}
    \left(
     \begin{array}{c}
       \hat{c}_{[1]}^{0}[n]_{1} \\
         \hat{c}_{[1]}^{1}[n]_{1}\\
     \end{array}
   \right)=\frac{1}{2} \tilde{\mathbf{M}}^{c}[-n]\cdot \left(
     \begin{array}{l}
      \hat{x}[n] \\
       \hat{x}[\vec{n}]
     \end{array}
   \right),\quad  \left(
     \begin{array}{l}
      \hat{x}[n] \\
       \hat{x}[\vec{n}]
     \end{array}
   \right)={\mathbf{M}}^{c}[n]\cdot \left(
     \begin{array}{c}
        \hat{c}_{[1]}^{0}[n]_{1} \\
         \hat{c}_{[1]}^{1}[n]_{1}\\
     \end{array}
   \right),
    \end{equation*}
    where $\vec{n}=n+{N}/{2}$.

Define the p-filters $\mathbf{q}^{l}_{\pm[1]}\srr \mathbf{h}^{j}_{[1]}\pm i\,\mathbf{g}^{j}_{[1]}=\psi^{p}_{[1],l}\pm i\,\varphi^{p}_{[1],l}={\Psi}^{p}_{\pm [1],l}, \;l=0,1.$
\ehh{qa_df} implies that their frequency response s are
\begin{eqnarray*}\label{Psi1_dfp}
  \hat{q}^{0}_{+[1]}[n]=\left\{
               \begin{array}{ll}
  (1+i)\sqrt{2}, & \hbox{if $n=0$;} \\
                  2\beta[n], & \hbox{if $0< n< N/2$;} \\
               0 & \hbox{if $ N/2\leq n<N$,}
               \end{array}
             \right.\quad
 \hat{q}^{1}_{+[1]}[n]=\left\{
               \begin{array}{ll}
-(1+i)\sqrt{2}, & \hbox{if  $n= N/2$;} \\
                   2\alpha[n], & \hbox{if $0< n< N/2$;} \\
               0, & \hbox{if $ N/2< n\leq N$.}
               \end{array}
             \right.\\\label{Psi1_dfm}
  \hat{q}^{0}_{-[1]}[n]=\left\{
               \begin{array}{ll}
  (1-i)\sqrt{2}, & \hbox{if $n=0$;} \\
                  2\beta[n], &  \hbox{if $ N/2<n<N$,} \\
               0 & \hbox{if $0< n\leq  N/2$;}
               \end{array}
             \right.\quad
 \hat{q}^{1}_{-[1]}[n]=\left\{
               \begin{array}{ll}
-(1-i)\sqrt{2}, & \hbox{if  $n= N/2$;} \\
                   2\alpha[n] & \hbox{if $ N/2< n\leq  N$;} \\
               0, & \hbox{if $0\leq n< N/2$.}
               \end{array}
             \right.
\end{eqnarray*}
Thus, the analysis  modulation matrices for the p-filters $\mathbf{q}^{l}_{\pm[1]}$ are
\begin{eqnarray}\label{aa_modma10p}
 \tilde{\mathbf{M}}_{+}^{q}[n]&=& \left(
                         \begin{array}{cc}
                            \hat{q}_{+[1]}^{0}[n] &   0\\
                             \hat{q}^{1}_{+[1]}[n] &   -\sqrt{2}(1+i)\,\delta[n-N/2]\\
                         \end{array}
                       \right)= \tilde{\mathbf{M}}[n]+i\, \tilde{\mathbf{M}}^{c}[n],\\\label{aa_modma10m}
 \tilde{\mathbf{M}}_{-}^{q}[n]&=& \left(
                         \begin{array}{cc}
                     (1-i)\sqrt{2}\delta[n]      &   \hat{q}_{-[1]}^{0}[n] \\
                             0 &  \hat{q}^{1}_{-[1]}[n]\\
                         \end{array}
                       \right)= \tilde{\mathbf{M}}[n]-i\, \tilde{\mathbf{M}}^{c}[n],
  \end{eqnarray}
   where the modulation matrix $\tilde{\mathbf{M}}[n]$ is  defined in Eq. \rf{sa_modma10T} and $\tilde{\mathbf{M}}^{c}[n]$ is  defined in \eh{aa_modma11}.
Application of the matrices $\tilde{\mathbf{M}}_{\pm}^{q}[n]$ to the vector $( \hat{x}[n] ,
       \hat{x}[\vec{n}])^{T}$ produces the vectors
\begin{equation}\label{mod_decAQ1}
    \left(
     \begin{array}{c}
       \hat{z}_{\pm[1]}^{0}[n]_{1} \\
         \hat{z}_{\pm[1]}^{1}[n]_{1}\\
     \end{array}
   \right)=\frac{1}{2}( \tilde{\mathbf{M}}_{\pm}^{q}[n])^{*}\cdot \left(
     \begin{array}{l}
      \hat{x}[n] \\
       \hat{x}[\vec{n}]
     \end{array}
   \right)= \left(
     \begin{array}{c}
       \hat{y}_{[1]}^{0}[n]_{1} \\
         \hat{y}_{[1]}^{1}[n]_{1}\\
     \end{array}
   \right)\mp i\, \left(
     \begin{array}{c}
       \hat{c}_{[1]}^{0}[n]_{1} \\
         \hat{c}_{[1]}^{1}[n]_{1}\\
     \end{array}
   \right).
 \end{equation}

\ehh{mod_decAQ1} implies that the inverse DFTs of the sequences $\hat{z}_{\pm[1]}^{\mu}[n]_{1},\;\mu=0,1,$ are
\begin{equation}\label{mod_decAQi}
  {z}_{\pm[1]}^{j}[l]=\left\langle \mathbf{x},\Psi^{p}_{\pm[1],j}[\cdot, -2l] \right\rangle=\sum_{k=0}^{N-1}x[k]\,\Psi^{p}_{\pm[1],j}[k -2l]^{*},\quad l=0,...,N/2-1.
\end{equation}

Define the matrices ${\mathbf{M}}_{\pm}^{q}[n]\srr\tilde{\mathbf{M}}_{\pm}^{q}[n]={\mathbf{M}}[n]\pm i\,{\mathbf{M}}^{c}[n]$ and apply these matrices to the vectors\\ $(\hat{z}_{\pm[1]}^{0}[n]_{1} ,
         \hat{z}_{\pm[1]}^{1}[n]_{1})^{T}$. Here   the modulation matrix  ${\mathbf{M}}[n]$ is defined in Eq.  \rf{sa_modma10T} and ${\mathbf{M}}^{c}[n]$ is  defined in \eh{aa_modma11}.
\bpp\label{pro:Mq_z}The following relations hold
\begin{eqnarray*}\label{Mq_z1}
 && {\mathbf{M}}_{\pm}^{q}[n]\cdot \left(
     \begin{array}{c}
        \hat{z}_{\pm[1]}^{0}[n]_{1} \\
         \hat{z}_{\pm[1]}^{1}[n]_{1}\\
     \end{array}
   \right)=\mathbf{M}[n]\cdot \left(
     \begin{array}{c}
        \hat{y}_{[1]}^{0}[n]_{1} \\
         \hat{y}_{[1]}^{1}[n]_{1}\\
     \end{array}
   \right) +{\mathbf{M}}^{c}[n]\cdot \left(
     \begin{array}{c}
        \hat{c}_{[1]}^{0}[n]_{1} \\
         \hat{c}_{[1]}^{1}[n]_{1}\\
     \end{array}
   \right) \\\nonumber&&\pm i \left( \mathbf{M}^{c}[n]\cdot \left(
     \begin{array}{c}
        \hat{y}_{[1]}^{0}[n]_{1} \\
         \hat{y}_{[1]}^{1}[n]_{1}\\
     \end{array}
   \right)   -\mathbf{M}[n]\cdot \left(
     \begin{array}{c}
        \hat{c}_{[1]}^{0}[n]_{1} \\
         \hat{c}_{[1]}^{1}[n]_{1}\\
     \end{array}
   \right)     \right)
    \\\label{Mq_z2}&&=2\left(\left(
     \begin{array}{l}
      \hat{x}[n] \\
       \hat{x}[n+N/2]
     \end{array}
   \right)
{\pm}i \,\left(
    \begin{array}{l}
      \hat{h}[n] \\
       \hat{h}[n+N/2]
     \end{array}
   \right)\right)
=2\left(
     \begin{array}{l}
      \hat{{x}}_{\pm}[n] \\
       \hat{{x}}_{\pm}[n+N/2]
     \end{array}
   \right)
,
\end{eqnarray*}
where $\mathbf{h}$ is the HT of the signal  $\mathbf{x}\in \Pi[N]$ and  $\mathbf{{x}}_{\pm}$ are the analytic  signals associated with  $\mathbf{x}$.
\epp

\proof In  Appendix.
\bdd\label{def:a_mvs} The matrices  $\tilde{\mathbf{M}}_{\pm}^{q}[n]$ and  ${\mathbf{M}}_{\pm}^{q}[n]$ are called the analysis  and synthesis  modulation matrices for the qWP transform, respectively.\edd
\br\label{rem:a_mvs}Successive application of the filter banks  $\tilde{\mathbf{H}}_{\pm}^{q}$ and  ${\mathbf{H}}_{\pm}^{q}$ defined by the analysis  and synthesis  modulation matrices  $\tilde{\mathbf{M}}_{\pm}^{q}[n]$ and  ${\mathbf{M}}_{\pm}^{q}[n]$, respectively, to a signal  $\mathbf{x}\in \Pi[N]$ produces the analytic  signals  $\mathbf{{x}}_{\pm}$ associated with  $\mathbf{x}$:
\begin{equation}\label{2MM}
 {\mathbf{H}}_{\pm}^{q}\cdot \tilde{\mathbf{H}}_{\pm}^{q}\cdot \mathbf{x}=2\bar{\mathbf{x}}_{\pm}\Longrightarrow \mathbf{x}=2\mathfrak{Re}( {\mathbf{H}}_{\pm}^{q}\cdot \tilde{\mathbf{H}}_{\pm}^{q}\cdot \mathbf{x}).
\end{equation}
\er
\bcc\label{cor:on_sys}A signal  $\mathbf{x}\in \Pi[N]$  is represented by the redundant   system
\begin{eqnarray*}\label{on_sys}
         x[k]&=&\frac{1}{2}\sum_{j=0}^{1}\sum_{l=0}^{N/2-1}\left(y_{[1]}^{j} [l]\psi_{[1],j}^{p}[k-2l] +
          c_{[1]}^{j} [l] \varphi_{[1],j}^{p}[k-2l]\right),\\\label{on_sys_yc}
          y_{[1]}^{j} [l]&=& \left\langle \mathbf{x},\, \psi_{[1],j}^{p}[\cdot-2l]   \right\rangle,\quad
         c_{[1]}^{j} [l]= \left\langle \mathbf{x},\, \varphi_{[1],j}^{p}[\cdot-2l]   \right\rangle.
         \end{eqnarray*}
         Thus, the system
        \begin{equation*}\label{on_sysTF}
\mathbf{F} \srr  \left\{ \left\{\psi_{[1],0}^{p}[\cdot-2l] \right\}\bigoplus  \left\{\psi_{[1],1}^{p}[\cdot-2l] \right\} \right\}\bigcup \left\{\left\{\varphi_{[1],0}^{p}[\cdot-2l] \right\}\bigoplus \left\{\varphi_{[1],1}^{p}[\cdot-2l] \right\} \right\},
\end{equation*}
whose components are orthonormal, form a tight frame of the space $ \Pi[N]$. Here $ l=0,..,N/2-1.$
\ecc
\subsection{Multi-level transforms}\label{sec:ss42}

%\subsubsection{Second-level transforms}\label{sec:sss421}
It was explained in Section \ref{sec:sss222} that the second-level transform  coefficients $\mathbf{y}_{[2]}^{\rho}$  are
\begin{eqnarray*}
 {y}_{[2]}^{\rho} [l]&=&   \sum_{n=0}^{N-1}x[n]\,{\psi}_{[2],\rho}^{p}[n-4l], \quad {\psi}_{[2],\rho}^{p}[n]  =\sum_{k=0}^{N/2-1}{h}_{[2]}^{\mu}[k] \, {\psi}_{[1],\lambda}^{p}[n-2k]\Longrightarrow\\
   {y}_{[2]}^{\rho} [l]&=&  \sum_{k=0}^{N/2-1}h_{[2]}^{\mu}[k-2l] \, y_{[1]}^{\lambda}[k], \quad \lambda,\mu=0,1,\;\rho=\left\{
                                                                                                                               \begin{array}{ll}
                                                                                                                                 \mu, & \hbox{if $\lambda=0$ ;} \\
                                                                                                                                 3-\mu, & \hbox{if $\lambda=1$.}
                                                                                                                               \end{array}
                                                                                                                             \right.
\end{eqnarray*}
The frequency response s of the p-filters  are $\hat{h}_{[2]}^{0}[n]=\beta[2n]$ and $\hat{h}_{[2]}^{1}[n]=\alpha[2n]$. The direct and inverse transforms
$\mathbf{y}_{[1]}^{\lambda}\longleftrightarrow\mathbf{y}_{[2]}^{2\lambda}\bigcup\mathbf{y}_{[2]}^{2\lambda+1}$ are  implemented using the analysis  and synthesis  modulation matrices $\tilde{\mathbf{M}}[2n]$ and $\mathbf{M}[2n]$, respectively.

The second-level transform  coefficients $\mathbf{c}_{[2]}^{\rho}$ are
\begin{eqnarray*}
 {c}_{[2]}^{\rho} [l]&=&   \sum_{n=0}^{N-1}x[n]\,{\varphi}_{[2],\rho}^{p}[n-4l], \quad {\varphi}_{[2],\rho}^{p}[n]  =\sum_{k=0}^{N/2-1}{h}_{[2]}^{\mu}[k] \, {\varphi}_{[1],\lambda}^{p}[n-2k]\Longrightarrow\\
   {c}_{[2]}^{\rho} [l]&=&  \sum_{k=0}^{N/2-1}h_{[2]}^{\mu}[k-2l] \, c_{[1]}^{\lambda}[k], \quad \lambda,\mu=0,1,\;\rho=\left\{
                                                                                                                               \begin{array}{ll}
                                                                                                                                 \mu, & \hbox{if $\lambda=0$ ;} \\
                                                                                                                                 3-\mu, & \hbox{if $\lambda=1$.}
                                                                                                                               \end{array}
                                                                                                                             \right.
\end{eqnarray*}
We emphasize that the p-filters $\mathbf{h}_{[2]}^{\mu}$ for the transform  $\mathbf{c}_{[1]}^{\lambda}\longleftrightarrow\mathbf{c}_{[2]}^{2\lambda}\bigcup\mathbf{c}_{[2]}^{2\lambda+1}$ are the same that the p-filters for the transform  $\mathbf{y}_{[1]}^{\lambda}\longleftrightarrow\mathbf{y}_{[2]}^{2\lambda}\bigcup\mathbf{y}_{[2]}^{2\lambda+1}$. Therefore, the direct and inverse transforms
$\mathbf{c}_{[1]}^{\lambda}\longleftrightarrow\mathbf{c}_{[2]}^{2\lambda}\bigcup\mathbf{c}_{[2]}^{2\lambda+1}$ are  implemented using the same analysis  and synthesis  modulation matrices $\tilde{\mathbf{M}}[2n]$ and $\mathbf{M}[2n]$. Apparently, it is the case also for the transforms $\mathbf{z}_{\pm[1]}^{\lambda}\longleftrightarrow\mathbf{z}_{\pm[2]}^{2\lambda}\bigcup\mathbf{z}_{\pm[2]}^{2\lambda+1}$. The transforms to subsequent decomposition   levels are implemented in an iterative way:

\begin{eqnarray*}\label{mod_decA2m}
    \left(
     \begin{array}{c}
       \hat{z}_{\pm[m+1]}^{\rho0}[n]_{m+1} \\
         \hat{z}_{\pm[m+1]}^{\rho1}[n]_{m+1}\\
     \end{array}
   \right)&=&\frac{1}{2} \tilde{\mathbf{M}}[-2^{m}n]\cdot \left(
     \begin{array}{l}
      \hat{z}_{\pm[m]}^{\lambda}[n]_{m} \\
        \hat{z}_{\pm[m]}^{\lambda}[\vec{n}]_{m}
     \end{array}
   \right),\\\nonumber\left(
     \begin{array}{l}
      \hat{z}_{\pm[m]}^{\lambda}[n]_{m} \\
        \hat{z}_{\pm[m]}^{\lambda}[\vec{n}]_{m}
     \end{array}
   \right)&=&{\mathbf{M}}[2^{m}n]\cdot\left(
     \begin{array}{c}
       \hat{z}_{\pm[m+1]}^{\rho0}[n]_{m+1} \\
         \hat{z}_{\pm[m+1]}^{\rho1}[n]_{m+1}\\
     \end{array}
   \right),
 \end{eqnarray*}
where
\(
\rho0=\left\{
              \begin{array}{ll}
                2\lambda, & \hbox{if $\lambda$ is even;} \\
                2\lambda+1, & \hbox{if $\lambda$ is odd,}
              \end{array}
            \right. \) and vice versa for $\rho1$, $\vec{n}=n+N/2^{m+1}$ and $m=1,...,M$.
By the application of the inverse DFT to the arrays $\left\{ \hat{z}_{\pm[m+1]}^{\rho}[n]_{m+1}\right\}$, we get the arrays\\ $\left\{ z_{\pm[m+1]}^{\rho}[k]=y_{[m+1]}^{\rho}[k]\pm i\,c_{[m+1]}^{\rho}[k]\right\}$ of the transform  coefficients with the qWPs $\Psi^{p}_{\pm[m+1],\rho}$.
\br\label{rem:zyc_cyz}By operating on  the transform   coefficients $\left\{ z_{\pm[m]}^{\rho}[k]\right\}$, we simultaneously operate on  the  arrays $\left\{ y_{[m]}^{\rho}[k]\right\}$ and $\left\{ c_{[m]}^{\rho}[k]\right\}$, which are the coefficients for the transforms with the DTSWPs $\psi^{p}_{[m],\rho}$ and cWPs $\varphi^{p}_{[m],\rho}$, respectively. The execution speed of the transform  with the qWPs  $\left\{ \Psi_{\pm[m]}^{p}\right\}= \psi_{[m]}^{p}\pm i \varphi_{[m]}^{p}$ is the same as the speed of the transforms with either WPs $\left\{ \psi_{[m]}^{p}\right\}$ or cWPs $\left\{ \varphi_{[m]}^{p}\right\}$.\er
 The transforms are executed in the spectral domain using the  FFT by the application of critically sampled two-channel filter banks to the half-band spectral components $(\hat{x}[n],\hat{x}[n+N/2])^{T}$  of a signal.

\begin{figure}%[H]
\centering
\includegraphics[width=3.in]{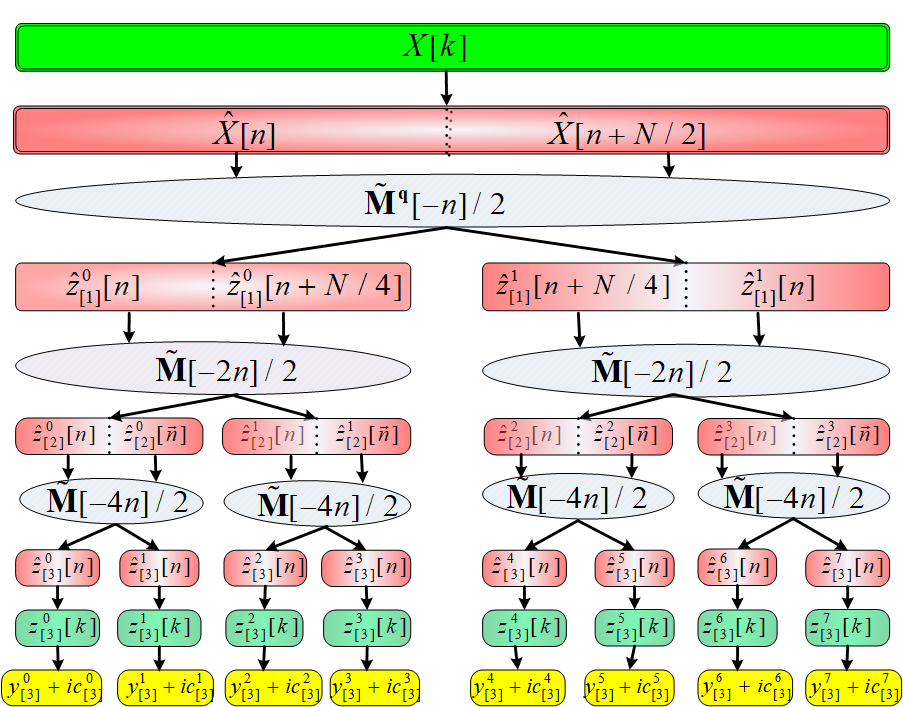}\hfill
\includegraphics[width=3.in]{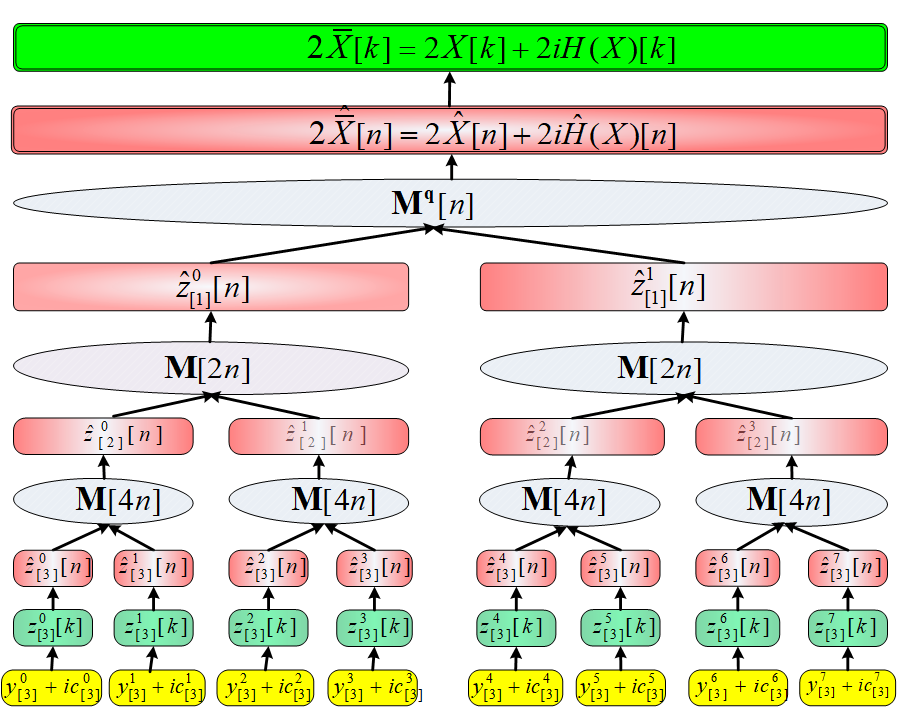}
\caption{ Left: Forward qWTP of a signal  $\mathbf{X}$ down to the third decomposition  level.  Right: Inverse transform   that results in restoration of the signal   $\mathbf{X}$ and its HT   $H(\mathbf{X})$. Here $\vec{n}$ means ${n}+N/8$ }
     \label{dia_wq_S}\label{dia_wq_A}
 \end{figure}

The diagrams in Fig.  \ref{dia_wq_S} illustrate the three-level forward and inverse qWP transforms of a signal  with quasi-analytic  wavelet packets, which use the analysis  $\tilde{\mathbf{M}}^{q}[n]$ and the synthesis  ${\mathbf{M}}^{q}[n]$ modulation matrices, respectively,  for the transforms to and from the first decomposition  level, respectively, and the modulation matrices  $\tilde{\mathbf{M}}[2^{m}n]$ and  ${\mathbf{M}}[2^{m}n]$ for the subsequent levels.

%\paragraph{Brief discussion}
\br\label{rem:recon}The decomposition  of a signal   $\mathbf{x}\in\Pi[N]$ down to the $M$-th level produces  $2MN$ transform  coefficients $\left\{ y_{[m]}^{\rho}[k]\right\}\bigcup\left\{ c_{[m]}^{\rho}[k]\right\}$. Such a redundancy provides many options for the signal  reconstruction. Some of them are listed below.
\begin{itemize}
  \item A basis compiled from either WPs $\left\{ \psi_{[m]}^{p}\right\}$ or $\left\{ \varphi_{[m]}^{p}\right\}$.
\begin{itemize}
  \item Wavelet basis.
  \item Best basis \cite{coiw1}, Local discriminant baseis \cite{sai}.
  \item WPs from a single decomposition  level.
\end{itemize}
  \item Combination of bases compiled from both $\left\{ \psi_{[m]}^{p}\right\}$ and $\left\{ \varphi_{[m]}^{p}\right\}$ WPs generates a tight frame of the space $\Pi[N]$  with redundancy rate 2. The bases for  $\left\{ \psi_{[m]}^{p}\right\}$ and $\left\{ \varphi_{[m]}^{p}\right\}$ can have a different structure.
  \item Frames with increased redundancy rate. For example, a combined reconstruction  from several decomposition  levels.
\end{itemize}
 \er
The collection of DTSWPs $\left\{ \psi_{[m]}^{p}\right\}$ and cWPs $\left\{ \varphi_{[m]}^{p}\right\}$, which originate from DTSs of different  orders $p$, provides a variety of waveforms that are (anti)symmetric, well localized in time domain. Any number of the  discrete  local vanishing moments can be achieved. The DFT spectra of the WPs are flat and the spectra shapes tend to rectangles when  the order $p$ increases. Therefore, they can be utilized as a collection of band-pass filters which produce a refined split of the frequency  domain into bands of different  widths. The (c)WPs can be used as testing waveforms for the signal  \aa, such as a dictionary for the Matching Pursuit procedures \cite{mal,azk_MP}.
\br\label{fpvm_rem}Since the magnitude spectra of the WPs $\psi_{[m],\lambda}^{p}$ and $ \varphi_{[m],\lambda}^{p}$ coincide, they have the same number of the  discrete  local vanishing moments.\er

\section{Two-dimensional complex wavelet packets}\label{sec:s5}
The 2D wavelet packets are defined by the tensor products of 1D WPs such that
\(%begin{equation*}\label{psipsi0}
  \psi_{[m],j ,l}^{p}[k,n]=\psi_{[m],j}^{p}[k]\,\psi_{[m], l}^{p}[n].
\)%end{equation*}
 The $2^{m}$-sample shifts of the DTSWPs $\left\{\psi_{[m],j ,l}^{p}\right\},\;j , l=0,...,2^{m}-1,$ in both directions form an orthonormal  basis for the space $\Pi[N,N]$ of arrays that are $N$-periodic  in both directions. The DFT spectrum  of such a WP is  concentrated in four symmetric spots in the frequency  domain.

Similar properties are inherent to the 2D cWPs such that
\(%begin{equation*}\label{phiphi}
  \varphi_{[m],j ,l}^{p}[k,n]=\varphi_{[m],j}^{p}[k]\,\varphi_{[m], l}^{p}[n].
\)%end{equation*}

\subsection{Design of  2D directional WPs}\label{sec:ss51}
\subsubsection{2D complex WPs and their spectra }\label{sec:sss511}
The DTSWPs  $\left\{\psi_{[m],j ,l}^{p}\right\}$ as well as the cWPs  $\left\{\varphi_{[m],j ,l}^{p}\right\}$  lack  the directionality property  which is needed in many applications that process 2D data.  However, real-valued 2D wavelet packets oriented in multiple directions  can be
derived from  tensor  products of complex  qWPs $\Psi_{\pm[m],\rho}^{p}$.

The complex 2D qWPs are defined  as follows:
%\begin{eqnarray*} \label{qwp_2d}
%\Psi_{++[m],j , l}^{p}[k,n] &\srr& \Psi_{+[m],j}^{p}[k]\,\Psi_{+[m], l}^{p}[n], \\\nonumber
%  \Psi_{+-[m],j ,l}^{p}[k,n] &\srr& \Psi_{+[m],j}^{p}[k]\,\Psi_{-[m], l}^{p}[n],
%\end{eqnarray*}
\[
\Psi_{++[m],j , l}^{p}[k,n] \srr \Psi_{+[m],j}^{p}[k]\,\Psi_{+[m], l}^{p}[n], \quad
  \Psi_{+-[m],j ,l}^{p}[k,n] \srr  \Psi_{+[m],j}^{p}[k]\,\Psi_{-[m], l}^{p}[n],
\]
where $  m=1,...,M,\;j ,l=0,...,2^{m}-1,$ and $k ,n=-N/2,...,N/2-1$.
The real  parts of these 2D qWPs are
\begin{equation}
\label{vt_pm}
\begin{array}{lll}
 \vartheta_{+[m],j ,l}^{p}[k,n] &\srr& \mathfrak{Re}(\Psi_{++[m],j ,l}^{p}[k,n]) =  \psi_{[m],j ,l}^{p}[k,n]-\varphi_{[m],j ,l}^{p}[k,n], \\
\vartheta_{-[m],j ,l}^{p}[k,n] &\srr&  \mathfrak{Re}(\Psi_{+-[m],j ,l}^{p}[k,n]) =  \psi_{[m],j ,l}^{p}[k,n]+\varphi_{[m],j ,l}^{p}[k,n],\\%\label{th_pm}
\end{array}
\end{equation}

\begin{SCfigure}
\centering
\caption{Magnitude spectra of 2D qWPs $\Psi_{++[2],j ,l}^{9}$ (left block of pictures)  and qWPs $\Psi_{+-[2],j ,l}^{9}$ (right block) from the second decomposition  level}
\includegraphics[width=1.9in]{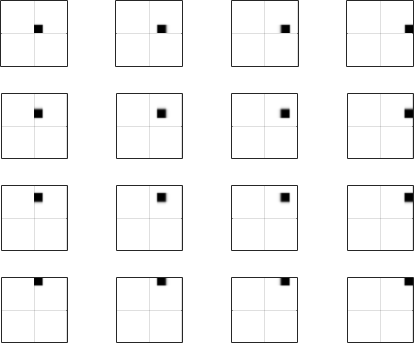}\quad\vline\quad
\includegraphics[width=1.9in]{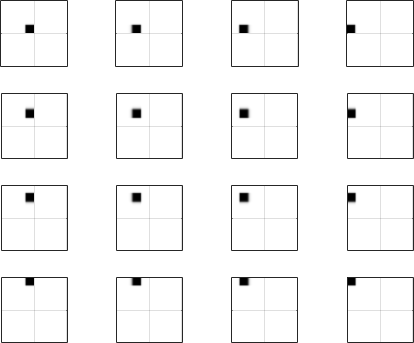}%
  \label{fppm_2}
\end{SCfigure}

The DFT spectra of the 2D qWPs $\Psi_{++[m],j ,l}^{p},\;j ,l=0,...,2^{m}-1,$ are the tensor products of the one-sided spectra of the qWPs:
%\begin{equation*}\label{fPsi_pp}
$\hat{ \Psi}_{++[m],j ,l}^{p}[p,q] =\hat{ \Psi}_{+[m],j}^{p}[p]\,\hat{\Psi}_{+[m], l}^{p}[q], $ 
   and, as such,  they fill the  quadrant $\mathbf{Q}_{0}$ of the frequency  domain, while the spectra of $\Psi_{+-[m],j ,l}^{p},\;j ,l=0,...,2^{m}-1,$ fill the  quadrant $\mathbf{Q}_{1}$ (see \eh{quadr}). Figure \ref{fppm_2} displays the magnitude spectra of the ninth-order 2D qWPs $\Psi_{++[2],j ,l}^{9}$ and $\Psi_{+-[2],j ,l}^{9}$ from the second decomposition  level, respectively.

\subsubsection{Directionality of real-valued 2D WPs}\label{sec:sss512}
It is seen in Fig. \ref{fppm_2} that the DFT spectra of the  qWPs $\Psi_{+\pm[m],j ,l}^{9}$ effectively occupy relatively small squares in the frequency  domain. For deeper decomposition  levels,  sizes of the corresponding  squares decrease on geometric progression. Such configurations of the spectra lead to the directionality of the real-valued 2D WPs $ \vartheta_{\pm[m],j ,l}^{p}$.

\begin{figure}
\centering
\includegraphics[width=2.5in]{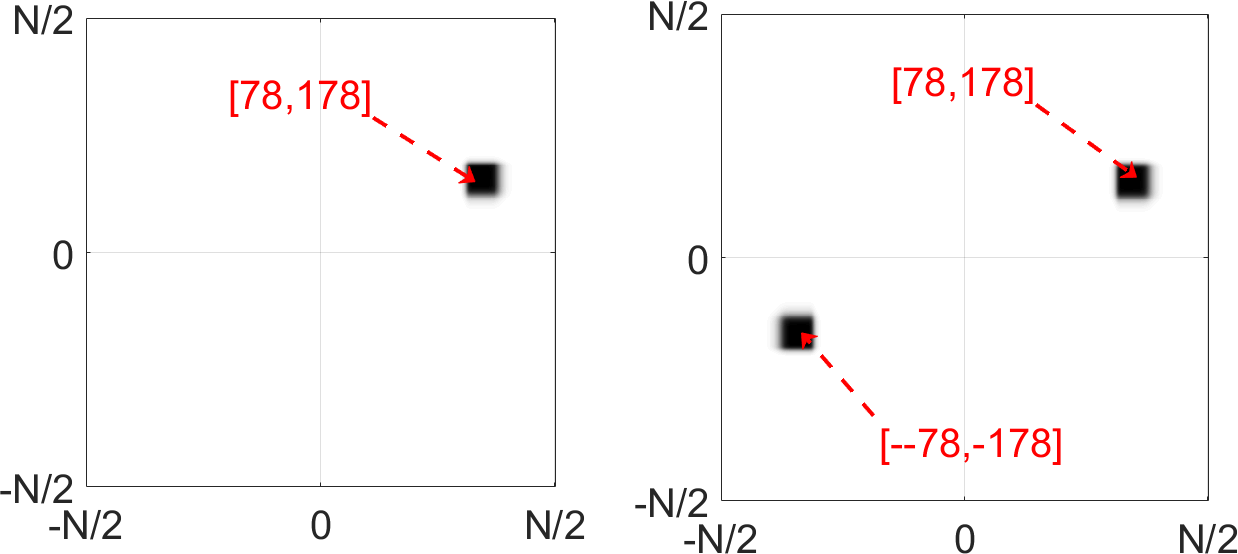}\quad\vline\quad
%\hfil
\includegraphics[width=3.6in]{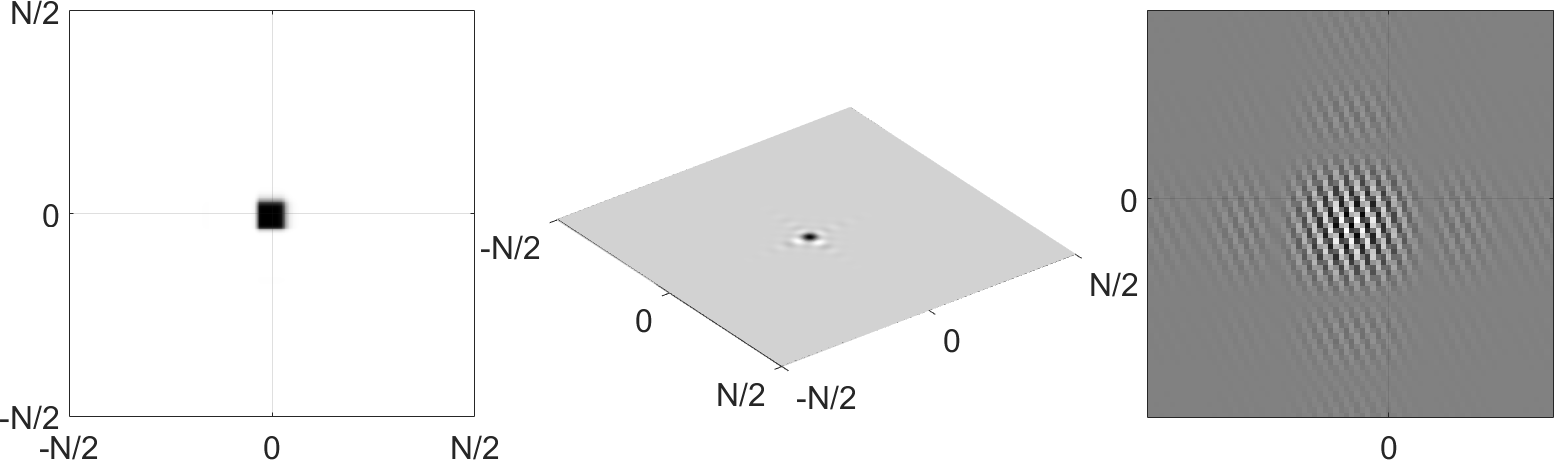}
\caption{Magnitude spectra of 2D qWP $\Psi[k,n]$ (left) and  $\mathfrak{Re}(\Psi)=\vartheta[k,n]$  (second from left). Center:  magnitude spectrum of  low-frequency  signal  $\underline{\vartheta}[k,n]$. Second from right: signal  $\underline{\vartheta}[k,n]$.  Right:  2D WP $\vartheta[k,n]$ (magnified)}
\label{78_178}
\end{figure}

Assume, for example, that $N=512, \;m=3,\; j=2,\;  l=5$ and denote $\Psi[k,n]\srr \Psi_{++[3],2 ,5}^{9}[k,n]$ and $ \vartheta[k,n]\srr\mathfrak{Re}(\Psi[k,n])$. The magnitude spectrum $\left|\hat{\Psi}[\kappa,\nu\right|$, displayed in Fig, \ref{78_178} (left), effectively occupies the square of size $40\times 40$ \emph{pixels} centered around the point $\mathbf{C}=[\kappa_{0},\nu_{0}]$, where $\kappa_{0}=78, \;\nu_{0}=178$. Thus, the WP $\Psi$ is represented by
\begin{eqnarray*}
\label{psi78_178}
\begin{array}{cc}
   \Psi[k,n] =& \frac{1}{N^{2}}\sum_{\kappa,\nu=0}^{N/2-1}\omega^{k\kappa+n\nu}\, \hat{\Psi}[\kappa,\nu]\approx{\omega^{\kappa_{0}k+\nu_{0}n}}\,\underline{\Psi}[k,n] , \\
   \underline{\Psi}[k,n] \srr & \frac{1}{N^{2}}\sum_{\kappa,\nu=-20}^{19}\omega^{k\kappa+n\nu} \, \hat{\Psi}[\kappa+\kappa_{0},\nu+\nu_{0}].
\end{array}
\end{eqnarray*}
Consequently, the real-valued  WP $\vartheta $, whose magnitude spectrum is displayed in Fig, \ref{78_178} (second from left), is represented as follows:
\begin{eqnarray*}
\label{th78_178}
  \vartheta[k,n]  \approx{\cos\frac{2\pi(\kappa_{0}k+\nu_{0}n)}{N}}\,\underline{\vartheta}[k,n] ,\quad \underline{\vartheta}[k,n] \srr\mathfrak{Re}(\underline{\Psi}[k,n]).
\end{eqnarray*}
The spectrum of the 2D signal  $\underline{\vartheta}$ comprises only  low frequencies in both directions and it does not have a directionality. But the 2D signal  $\cos\frac{2\pi(\kappa_{0}k+\nu_{0}n)}{N}$ is oscillating in the direction of the vector $\vec{V}_{++[2],2 ,5}=178\vec{i}+78\vec{j}$. The 2D WP $\vartheta[k,n]$ is well localized in the spatial domain as is seen from \eh{vt_pm} and the same is true for the low-frequency  signal  $\underline{\vartheta}$. Therefore, WP $\vartheta[k,n]$  can be regarded as the directional cosine modulated  by the localized low-frequency  signal  $\underline{\vartheta}$.

The same arguments are applicable to the 2D WPs $\vartheta_{-[m],j ,l}^{p}[k,n] = \mathfrak{Re}(\Psi_{+-[m],j ,l}^{p}[k,n]) $. Figure \ref{78_178} displays the low-frequency  signal  $\underline{\vartheta}$, its magnitude spectrum and the  2D WP $\vartheta[k,n]$.

\begin{SCfigure}%[H]
\centering
\caption{WPs $\vartheta_{+[2],j ,l}^{9}$ from the second decomposition  level and their magnitude spectra}
\includegraphics[width=2.3in]{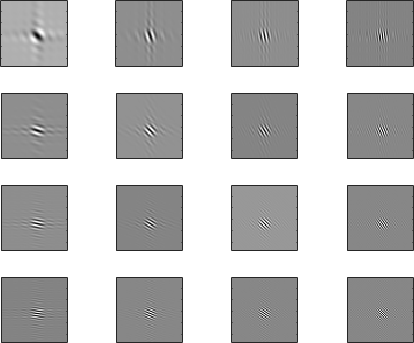}
\quad
\includegraphics[width=2.3in]{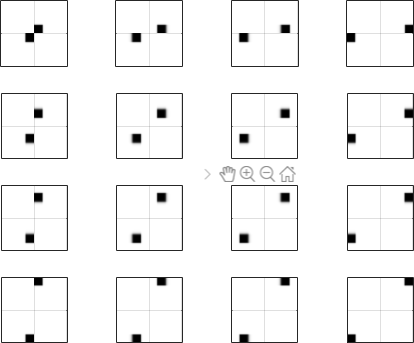}%
\
\label{pp_2_2d}
\end{SCfigure}

Figure \ref{pp_2_2d} displays WPs  $\vartheta_{+[2],j ,l}^{9},\;j,l=0,1,2,3,$ from the second decomposition  level and their magnitude spectra.

\begin{SCfigure}%[H]
\centering
\caption{WPs $\vartheta_{-[2],j ,l}^{9}$ from the second decomposition  level and their magnitude spectra}
\includegraphics[width=2.3in]{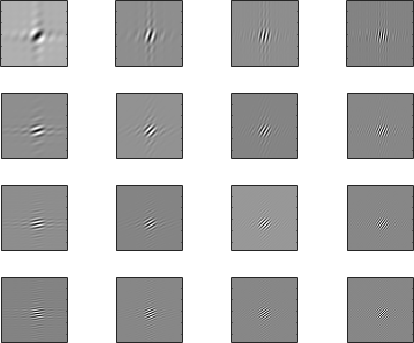}
\quad
\includegraphics[width=2.3in]{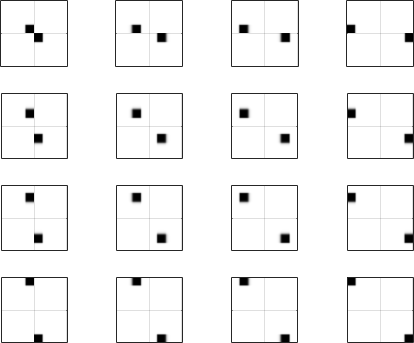}%
\label{pm_2_2d}
\end{SCfigure}

Figure \ref{pm_2_2d} displays WPs  $\vartheta_{-[2],j ,l}^{9},\;j,l=0,1,2,3,$ from the second decomposition  level and their magnitude spectra.

\begin{SCfigure}%[H]
\centering
\caption{WPs $\vartheta_{+[3],j ,l}^{9}$  (left) and $\vartheta_{-[3],j ,l}^{9}$  (right) from the third decomposition  level }
\includegraphics[width=2.6in]{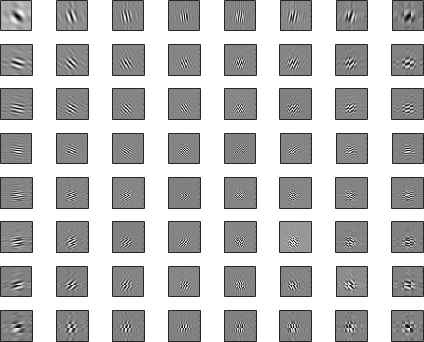}\quad\vline\quad
%\hfil
\includegraphics[width=2.6in]{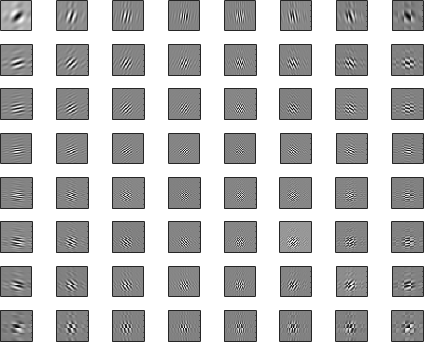}%
\label{pp_3_2d}
\end{SCfigure}

\br\label{direc_rem}Note that orientations of the vectors $\vec{V}_{++[m],j, l}$ and $\vec{V}_{++[m],j+1, l+1}$ are approximately the same. These vectors   determine the orientations of the WPs $\vartheta_{+[m],j ,l}^{p}$ and $\vartheta_{+[m],j+1 , l+1}^{p}$, respectively. Thus, these WPs have approximately the same orientation although they differ by the oscillation frequencies.  Consequently, the WPs from the $m$-th decomposition  level are oriented in $2^{m+1}-1$ different  directions. The same is true for the WPs $\vartheta_{-[m],j ,l}^{p}$. Thus, altogether, at the level $m$ we have WPs oriented in $2(2^{m+1}-1)$ different  directions.  It is seen in Figs. \ref{pp_2_2d}, \ref{pm_2_2d} and in Fig. \ref{pp_3_2d}, which displays the WPs $\vartheta_{\pm[3],j ,l}^{p}$.\er

\section{Implementation of 2D qWP transforms}\label{sec:s6}
The  spectra of the real-valued 2D WPs $\left\{\vartheta_{+[m],j ,l}^{p}\right\},\;j ,l=0,...,2^{m}-1$, and  $\left\{\vartheta_{-[m],j ,l}^{p}\right\}$ fill the pairs  of quadrant  $\mathbf{Q}_{+}\srr\mathbf{Q}_{0}\bigcup\mathbf{Q}_{2}$ and $\mathbf{Q}_{-}\srr\mathbf{Q}_{1}\bigcup\mathbf{Q}_{3}$ (see \eh{quadr}), respectively (Figs. \ref{pp_2_2d} and \ref{pm_2_2d}).

By this reason, none linear combination of the WPs $\left\{\vartheta_{+[m],j ,l}^{p}\right\}$  and their shifts can serve as a basis for the signal  space $\Pi[N,N]$. The same is   true for WPs $\left\{\vartheta_{-[m],j ,l}^{p}\right\}$. However,  combinations of  the WPs $\left\{\vartheta_{\pm[m],j ,l}^{p}\right\}$ provide frames of the  space $\Pi[N,N]$.
\subsection{One-level 2D transforms}\label{sec:ss61}
The one-level 2D qWP transforms of a signal  $\mathbf{X}=\left\{X[k,n] \right\}\in\Pi[N,N]$ are implemented by a  tensor-product scheme.
Denote by  $H^{v}(\mathbf{S})$ and $H^{h}(\mathbf{S})$ results of application of the Hilbert transforms to columns and rows of  a  2D signal    $\mathbf{S}$, respectively.

%\subsubsection{Direct transforms with qWPs $\Psi^{p}_{+\pm[1]}$}\label{sec:sss611}
  Denote by $\tilde{\mathbf{T}}_{\pm}^{h}$ the 1D transforms of row signals from $\Pi[N]$ with the analysis  modulation matrices $\tilde{\mathbf{M}}_{\pm}^{q}$ which are defined in \eh{aa_modma10p}. Application of these transforms  to rows of a signal   \textbf{X} produces the coefficient  arrays
\begin{eqnarray*}\label{tTh+x}
        \tilde{\mathbf{T}}_{+}^{h}\mathbf{\cdot}\mathbf{X} &=&  \left(\zeta_{+}^{0},\zeta_{+}^{1}\right),\quad \zeta_{+}^{j}[k,n]=\eta^{j}[k,n]-i\,\xi^{j}[k,n],
                                                             \\\nonumber \tilde{\mathbf{T}}_{-}^{h}\mathbf{\cdot}\mathbf{X} &=&  \left(\zeta_{-}^{0},\zeta_{-}^{1}\right),\quad \zeta_{-}^{j}[k,n]=\eta^{j}[k,n]+i\,\xi^{j}[k,n]=(\zeta_{+}^{j}[k,n])^{*},
                                                             \\\nonumber
        {\eta}^{j}[k,n] &=&\left\langle \mathbf{X}[k,\cdot],{\psi}^{p}_{[1],j}[\cdot -2n]\right\rangle,\quad  {\xi}^{j}[k,n]=\left\langle \mathbf{X}[k,\cdot],{\varphi}^{p}_{[1],j}[\cdot -2n]\right\rangle,\;j=0,1.
      \end{eqnarray*}
Here $\eta^{j}$ and  $\xi^{j}$ are real-valued arrays of size $N\times N/2$. Obviously we have
\begin{equation}\label{htvx}
 \tilde{\mathbf{T}}_{\pm}^{h}\mathbf{\cdot}H^{v}(\mathbf{X})=  \left(H^{v}(\zeta_{\pm}^{0}),H^{v}(\zeta_{\pm}^{1})\right), \quad H^{v}(\zeta_{\pm}^{j})= H^{v}(\eta^{j})\mp H^{v}(\xi^{j}).
\end{equation}

Denote by $\tilde{\mathbf{T}}_{+}^{v}$     the direct  1D transform  determined by the modulation matrix $\tilde{\mathbf{M}}_{+}^{q}$  applicable to columns of the corresponding signals. The next step of  the tensor product transform  consists of  the application of the  1D transform  $\tilde{\mathbf{T}}_{+}^{v}$ to  columns of the arrays  ${\zeta}^{j},\;j=0,1.$
\begin{eqnarray*}
  \tilde{\mathbf{T}}_{+}^{v}\mathbf{\cdot}\zeta_{+}^{j} =  \tilde{\mathbf{T}}_{+}^{v}\mathbf{\cdot}\,\eta^{j}- i \tilde{\mathbf{T}}_{+}^{v}\mathbf{\cdot}\,\xi^{j} =\mathbf{Z}_{+[1]}^{j}, \quad
   \tilde{\mathbf{T}}_{+}^{v}\mathbf{\cdot}\zeta_{-}^{j}=\tilde{\mathbf{T}}_{+}^{v}\mathbf{\cdot}\,\eta^{j}+ i \tilde{\mathbf{T}}_{+}^{v}\mathbf{\cdot}\,\xi^{j} =\mathbf{Z}_{-[1]}^{j}.
\end{eqnarray*}

Denote by $\mathbf{T}_{+}^{v}$ the 1D inverse transform  with the synthesis  modulation matrix ${\mathbf{M}}_{+}^{q}$ applicable to columns of the coefficient  arrays.
\begin{eqnarray*}
  \mathbf{T}_{+}^{v}\mathbf{\cdot}\mathbf{Z}_{+[1]}^{j} &=& 2(\eta^{j}+iH^{v}(\eta^{j}))-2i(\xi^{j}+iH^{v}(\xi^{j}))=2(\zeta_{+}^{j}+ iH^{v}(\zeta_{+}^{j})),\\
 \mathbf{T}_{+}^{v}\mathbf{\cdot}\mathbf{Z}_{-[1]}^{j} &=& 2(\eta^{j}+iH(\eta^{j})+2i(\xi^{j}+iH(\xi^{j}))=2(\zeta_{-}^{j}+ iH^{v}(\zeta_{-}^{j})).
\end{eqnarray*}
Denote by ${\mathbf{T}}_{\pm}^{h}$ the 1D inverse transforms with the synthesis  modulation matrices ${\mathbf{M}}_{\pm}^{q}$. %Due to Proposition \ref{pro:Mq_z},
Application of these transforms  to rows of the coefficient  arrays  $\zeta_{\pm}=\left(\zeta_{\pm}^{0},\zeta_{\pm}^{1}\right)$, respectively, produces the  2D  analytic  signals: $
{\mathbf{T}}_{\pm}^{h}\mathbf{\cdot}(\zeta_{\pm}^{0},\zeta_{\pm}^{1})=2(\mathbf{X}\pm i\,H^{h}(\mathbf{X}))$.

\ehh{htvx} implies that application of the transforms  ${\mathbf{T}}_{\pm}^{h}$ to rows of the  arrays  $H^{v}(\zeta_{\pm})\srr\ \left(H^{v}(\zeta_{\pm}^{0}),H^{v}(\zeta_{\pm}^{1})\right)$, respectively, produces the  2D  analytic  signals: $
{\mathbf{T}}_{\pm}^{h}\mathbf{\cdot}\left(H^{v}(\zeta_{\pm}^{0}),H^{v}(\zeta_{\pm}^{1})\right))=2(\mathbf{G}\pm i\,H^{h}(\mathbf{G}))$, where $\mathbf{G}=H^{v}(\mathbf{X}).$
Consequently,
 \begin{eqnarray}\nonumber
  \mathbf{X}_{+}&\srr&{\mathbf{T}}_{+}^{h}\mathbf{\cdot}  \mathbf{T}_{+}^{v}\mathbf{\cdot}\mathbf{Z}_{+[1]}^{j} =
4\left(\mathbf{X}+i\,H^{h}(\mathbf{X})+i\mathbf{G}-H^{h}(\mathbf{G})\right),\\\nonumber
\mathbf{X}_{-}&\srr&{\mathbf{T}}_{-}^{h}\mathbf{\cdot}  \mathbf{T}_{+}^{v}\mathbf{\cdot}\mathbf{Z}_{-[1]}^{j}
=4\left(\mathbf{X}-i\,H(\mathbf{X})+i\mathbf{G}+H^{h}(\mathbf{G})\right).\\\label{Th_zahza}
\mathbf{X}&=&\mathfrak{Re}\left(\frac{ \mathbf{X}_{+}+ \mathbf{X}_{-}}{8}\right).
 \end{eqnarray}

\begin{SCfigure}%[H]
\centering
\caption{Left to right: 1. Image  $\mathfrak{Re}(\mathbf{X}_{+})$. 2. Its magnitude DFT spectrum. 3.Image  $\mathfrak{Re}(\mathbf{X}_{-})$. 4. Its magnitude DFT spectrum }
\includegraphics[width=4.2in]{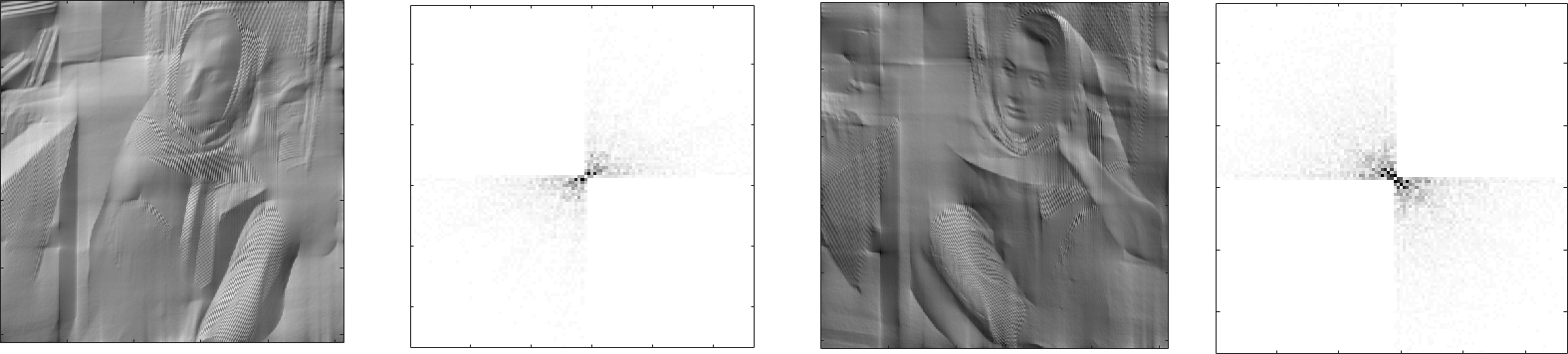}
\label{xp_xm_x5T2}
\end{SCfigure}

Figure \ref{xp_xm_x5T2} illustrates the  image ``Barbara" restoration by the 2D signals $\mathfrak{Re}(\mathbf{X}_{\pm})$. The signal  $\mathfrak{Re}(\mathbf{X}_{+})$ captures edges oriented to \emph{north-east}, while $\mathfrak{Re}(\mathbf{X}_{-})$ captures edges oriented to \emph{north-west}. The signal  $\tilde{\mathbf{X}}=\mathfrak{Re}(\mathbf{X}_{+}+\mathbf{X}_{-})/8$ perfectly restores the image achieving  PSNR=313.8596 dB.

\subsection{Multi-level 2D transforms}\label{sec:ss62}
It was established in Section \ref{sec:ss42} that the 1D qWP transforms of a signal  $\mathbf{x}\in\Pi[N]$  to the second and further  decomposition  levels are implemented by the successive  application of
the filter banks,  that  are determined by their analysis  modulation matrices $\tilde{\mathbf{M}}[2^{m}n],\;m=1,...,M-1,$ to the coefficient  arrays $\mathbf{z}_{\pm[m]}^{\lambda}$. The transforms  applied to the arrays $\mathbf{z}_{\pm[m]}^{\lambda}$  produce the arrays $\mathbf{z}_{\pm[m+1]}^{\rho}$, respectively. The inverse transform  consists of  the iterated  application of the
filter banks that are determined by their synthesis  modulation matrices ${\mathbf{M}}[2^{m}n],\;m=1,...,M-1,$ to the coefficient  arrays $\mathbf{z}_{\pm[m+1]}^{\rho}$.  In that way the first-level coefficient  arrays $\mathbf{z}_{\pm[1]}^{\lambda},\;\lambda=0,1$ are restored\footnote{The matrices $\tilde{\mathbf{M}}[n]$ and ${\mathbf{M}}[n]$ are defined in \eh{sa_modma10T}.}.

The tensor-product  2D transform  of a signal  $\mathbf{X}\in\Pi[N,N]$ consists of the subsequent application of the 1D transforms to columns and rows of the signal  and coefficient  arrays. By application of filter banks, which are determined by the analysis  modulation matrix $\tilde{\mathbf{M}}[2n]$ to columns and rows of  the coefficient  arrays $\mathbf{Z}_{\pm[1]}^{j,l}$, we derive four second-level arrays
$\mathbf{Z}_{\pm[2]}^{\rho,\tau},\;\rho=2j,2j+1;\;\tau=2l,2l+1$. The arrays $\mathbf{Z}_{\pm[1]}^{j,l}$ are restored by the application of the
filter banks that are determined by their synthesis  modulation matrices ${\mathbf{M}}[2n] $ to rows and columns of the coefficient  arrays $\mathbf{Z}_{\pm[2]}^{\rho,\tau},\;\rho=2j,2j+1;\;\tau=2l,2l+1$. The transition from the second to further levels and back are executed similarly using the modulation matrices $\tilde{\mathbf{M}}[2^{m}n]$ and ${\mathbf{M}}[2^{m}n]$, respectively. The inverse transforms produce the coefficient  arrays $\mathbf{Z}_{\pm[1]}^{j,l},\;j,l=0,1,$ from which the signal  $\mathbf{X}\in\Pi[N,N]$ is restored using the synthesis  modulation matrices ${\mathbf{M}}_{\pm}^{q}[n] $ as it is explained in Section \ref{sec:ss61}.

All the computations are implemented in the frequency  domain using the FFT.
%\br\label{rem:dob_tre}
\paragraph{Summary}
The 2D qWP processing  of a signal  $\mathbf{X}\in\Pi[N,N]$ is implemented by a dual-tree scheme. The first step produces two sets of the coefficient  arrays: $\mathbf{Z}_{+[1]}=\left\{\mathbf{Z}_{+[1]}^{j,l}\right\}, \;j,l,=0,1,$ which are derived using the analysis  modulation matrix $\tilde{\mathbf{M}}_{+}^{q}[n]$ for the row and column transforms, and $\mathbf{Z}_{-[1]}=\left\{\mathbf{Z}_{-[1]}^{j,l}\right\}, \;j,l,=0,1,$ which are derived using the analysis  modulation matrices  $\tilde{\mathbf{M}}_{+}^{q}[n]$ for the  column and $\tilde{\mathbf{M}}_{-}^{q}[n]$ for the row transforms.  Further decomposition  steps are implemented in parallel on the sets $\mathbf{Z}_{+[1]}$ and  $\mathbf{Z}_{-[1]}$  using the same analysis  modulation matrices $\tilde{\mathbf{M}}[2^{m}n]$, thus producing two multi-level sets of the coefficient  arrays $\left\{\mathbf{Z}_{+[m]}^{j,l}\right\}$ and $\left\{\mathbf{Z}_{-[m]}^{j,l}\right\},\;m=2,...,M,\;j,l=0,2^{m}-1$.

By parallel implementation of the inverse transforms on the coefficients from the   sets  $\left\{\mathbf{Z}_{+[m]}^{j,l}\right\}$ and $\left\{\mathbf{Z}_{-[m]}^{j,l}\right\}$   using the same synthesis  modulation matrix ${\mathbf{M}}[2^{m}n]$, the sets  $\mathbf{Z}_{+[1]}$ and  $\mathbf{Z}_{-[1]}$ are restored, which, in turn, provide the signals   $\mathbf{X}_{+}$ and  $\mathbf{X}_{-}$, using the synthesis  modulation matrices  ${\mathbf{M}}_{+}^{q}[n]$ and ${\mathbf{M}}_{-}^{q}[n]$, respectively. Typical  signals   $\mathfrak{Re}(\mathbf{X}_{\pm})$ and their DFT spectra are displayed in Fig. \ref{xp_xm_x5T2}.

Prior to the reconstruction, some structures, possibly different, are defined in the sets  $\left\{\mathbf{Z}_{+[m]}^{j,l}\right\}$ and $\left\{\mathbf{Z}_{-[m]}^{j,l}\right\},\;m=1,...M,$ (for example, 2D wavelet   or Best Basis structures) and some manipulations on the coefficients, (for example, thresholding, $l_1$ minimization) are executed.

\section{Discussion}\label{sec:s7} The paper describes the design of one- and two-dimensional quasi-analytic  WPs (qWPs) originating from polynomial  splines of arbitrary order and corresponding transforms. The qWP transforms operate in spaces of periodic  signals. Seemingly, the requirement of periodicity imposes some limitations on the scope of signals available for processing, but actually these limitations are easily circumvented by symmetrical extension of images  beyond the boundaries before processing  and shrinkage  to the original size after that. On the other hand, the periodic  setting provides a lot of substantial opportunities for the design and implementation of   WP transforms. The 2D qWPs possess the following properties:
\begin{description}
\item[-] The qWP transforms  provide a variety of 2D waveforms oriented in multiple directions. For example, fourth-level qWPs are oriented in 62 different  directions.
\item[-] The waveforms  are close to directional cosines with a variety of frequencies  modulated by spatially  localized low-frequency 2D signals and can have any number of local vanishing moments.
\item[-] The DFT spectra of the waveforms produce a refined tiling of the frequency  domain.
\item[-] Fast implementation of the transforms by using the FFT enables us to use the transforms with increased redundancy.
\end{description}

\begin{figure}%[H]
%\centering
\includegraphics[width=4.0in]{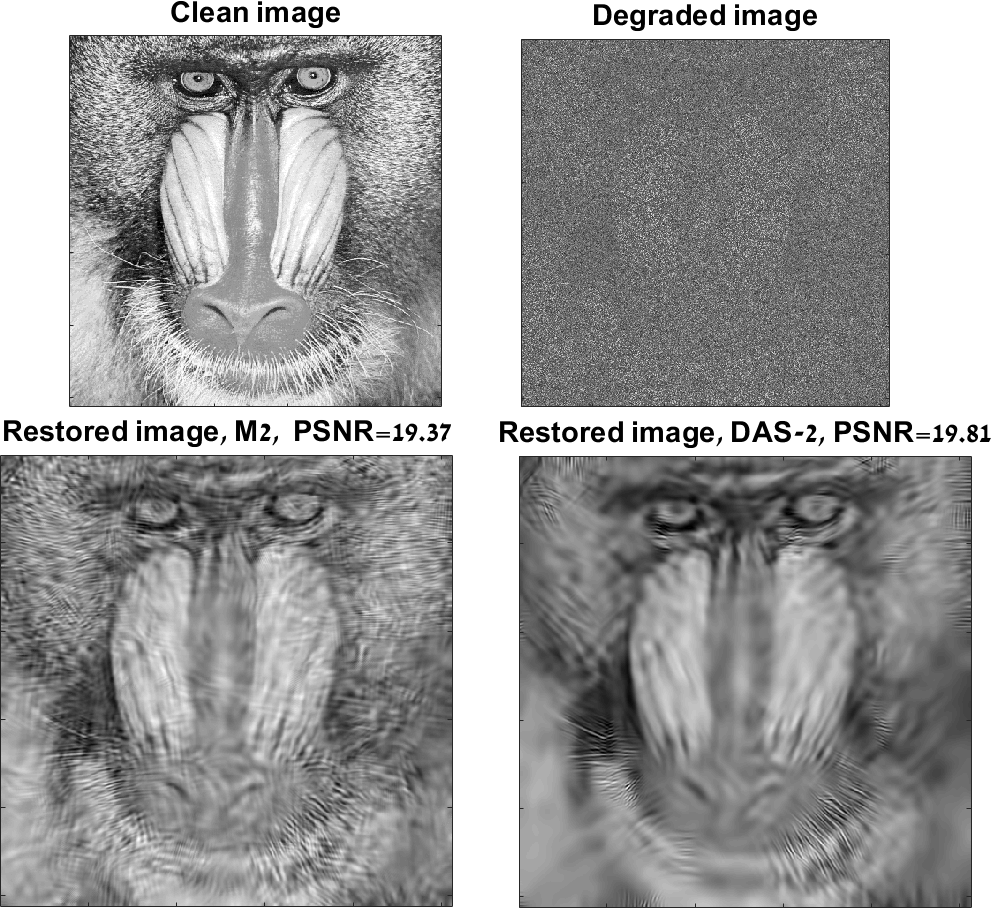}\hfill
\includegraphics[width=1.8in]{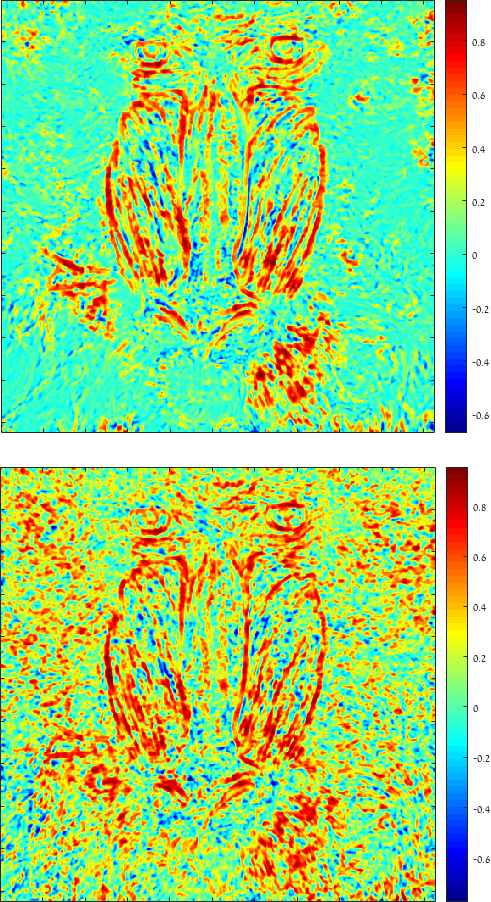}
\caption{Restoration of ``Mandrill" image. Left: Top left: clean image. Top right: image degraded by missing  80\% of its pixels   with additive    Gaussian noise with $\sigma=50$ dB.
Bottom left: \textbf{M2} restoration, PSNR=19.37 dB. Bottom right: DAS-2 restoration, PSNR=19.81 dB. Right:  SSIM map of images restored by  DAS-2 (top ), SSIM=0.1414, and by \textbf{M2} (bottom), SSIM=0.2185  }
\label{mand4_50}
\end{figure}

 The above listed properties of qWP transforms proved to be indispensable while dealing with  image processing   problems. Multiple experiments on image  denoising  and inpainting, whose results will be reported in our forthcoming publications,  demonstrate that qWP-based methods are quite competitive with the best  state-of-the-art   algorithms.
 Due to a variety of orientations, the qWPs  capture edges even in severely degraded images and their oscillating structures with  a variety of frequencies enable to  recover thin structures.  This fact is  illustrated in Fig. \ref{mand4_50}, which displays the restoration result  of the `Mandrill" image from the input where 80\% of its pixels are missing and additive Gaussian noise with $\sigma=50$ dB is present.  The result is compared with the output from  DAS-2 algorithm  (\cite{che_zhuang}). The output from DAS-2 has PSNR=19.81 dB compared to 19.37  dB produced by the qWP-based method designated by \textbf{M2}. However,  the Structural Similarity Index (SSIM) for the \textbf{M2}  restoration is 0.2185 compared to 0.1414 for DAS-2. The SSIM maps for \textbf{M2} and DAS-2  significantly differ from each other.
Figure \ref{mand4_50} is a good illustration to the fact that the SSIM  has much more informative characteristics than what PSNR provides.

Summarizing, by having such a versatile and
flexible tool at hand, %and being encouraged by the results achieved in solving the denoising  and inpainting problems,
we are in a position to address
multiple data processing problems such as image deblurring,  superresolution, segmentation and classification and target detection (here the directionality is of utmost importance). The  3D directional wavelet packets, whose design  is  underway, may be beneficial for seismic and  hyper-spectral processing.

We did not compare the qWP-based methods performance with the performance  of the schemes based on the deep learning (DL). However, we believe that  the designed directional qWPs can boost image processing   methods that are based on the Deep Learning by serving as a powerful tool for extraction of characteristic features  from images. This will be explored in
 our future work.
%\bibliographystyle{plain}
%   % mathematics and physical sciences
%\bibliography{BookBib_TBS}
%\end{document}
\include{ANA_TAPA}
\paragraph{Acknowledgment}
This research was partially supported by the Israel Science Foundation (ISF, 1556/17),
Supported by Len Blavatnik and the Blavatnik Family Foundation,
Israel Ministry of Science Technology and Space 3-16414, 3-14481 and by Academy of Finland (grant 311514).

%%%%%%%%%%%%%

\bibliographystyle{plain}
   % mathematics and physical sciences
\bibliography{BookBib_TBSA}
\end{document}